\documentclass{agtart_a}
\pdfoutput=1

\usepackage{graphicx}
\usepackage[all]{xy}

\title{A volume form on the $\mathrm{SU}(2)$--representation\\space of knot groups}
\author{J\'er\^ome Dubois}
\givenname{J\'er\^ome}
\surname{Dubois}
\address{Section de Math\'ematiques\\ 
Universit\'e de Gen\`eve CP 64\\\newline
2--4 Rue du Li\`evre\\ 
CH-1211 Gen\`eve 4\\
Switzerland}
\email{Jerome.Dubois@math.unige.ch}



\received{24 September 2004}
\revised{25 August 2005}
\accepted{27 December 2005}
\proposed{}
\seconded{}
\publishedonline{12 March 2006}
\published{12 March 2006}

\volumenumber{6}
\issuenumber{}
\publicationyear{2006}
\papernumber{13}
\startpage{373}
\endpage{404}

\doi{}
\MR{}
\Zbl{}

\arxivreference{math.GT/0409529}

\let\xysavmatrix\xymatrix
\def\xymatrix{\disablesubscriptcorrection\xysavmatrix}
\AtBeginDocument{\let\bar\wbar\let\tilde\wtilde\let\hat\what}

\newcommand{\ie}{ie, }
\newcommand{\cf}{cf~}

\newcommand{\SU}{\mathrm{SU}(2)}

\newcommand{\SO}{\mathrm{SO}(3)}
\newcommand{\su}{\mathfrak{su}(2)}

\newcommand{\ii}{\mathbf{i}}
\newcommand{\jj}{\mathbf{j}}

\newcommand{\kk}{\mathbf{k}}
\newcommand{\I}{\mathbf{1}}
\newcommand{\ZZ}{\mathbb{Z}}
\newcommand{\CC}{\mathbb{C}}
\newcommand{\IR}{\mathbb{R}}

\newcommand{\Reg}[1]{\mathcal{R}\mathrm{eg}(#1)}
\newcommand{\tangent}[2]{T_{#1} #2}

\newcommand{\vol}[1]{\omega^{#1}}
\newcommand{\fox}[2]{\frac{\partial #1}{\partial #2}}

\newcommand{\diff}[3]{\left. \frac{d #1}{d #2}\right\vert_{#3}}

\newcommand{\bord}{\partial}
\newcommand{\lk}{\ell\mathit{k}}

\newcommand{\rrk}{\mathop{\mathrm{rk}}\nolimits}

\makeatletter
\def\cnewtheorem#1[#2]#3{\newtheorem{#1}{#3}[section]
\expandafter\let\csname c@#1\endcsname\c@theorem}

\theoremstyle{plain}
\newtheorem{theorem}{Theorem}[section]
\newtheorem*{theorem*}{Theorem}

\newtheorem*{maintheorem*}{Main Theorem}
\cnewtheorem{prop}[theorem]{Proposition}%[chapter]
\cnewtheorem{corollary}[theorem]{Corollary}%[chapter]
\cnewtheorem{lemma}[theorem]{Lemma}%[chapter]
\cnewtheorem{claim}[theorem]{Claim}%[chapter]
\cnewtheorem{fact}[theorem]{Fact}
\theoremstyle{definition}
\newtheorem{example}{Example}%[chapter]
\newtheorem*{example*}{Example}
\newtheorem*{notation*}{Notation}
\newtheorem*{convention}{Convention}
\theoremstyle{remark}

\makeatother  %  move after \newtheorem block

\begin{document}

\begin{asciiabstract}
For a knot K in S^3 we construct according to Casson--or more
precisely taking into account Lin and Heusener's further works--a
volume form on the SU(2)-representation space of the group of
K.  We prove that this volume form is a topological knot invariant
and explore some of its properties.
\end{asciiabstract}

\begin{htmlabstract}
For a knot K in <i>S&sup3;</i> we construct according to Casson &ndash;
or more precisely taking into account Lin and Heusener's further works
&ndash; a volume form on the SU(2)&ndash;representation space of the group
of K.  We prove that this volume form is a topological knot invariant
and explore some of its properties.
\end{htmlabstract}

\begin{webabstract} 
For a knot $K$ in $S^3$ we construct according to Casson---or more
precisely taking into account Lin and Heusener's further works---a
volume form on the $\mathrm{SU}(2)$--representation space of the group of
$K$.  We prove that this volume form is a topological knot invariant
and explore some of its properties.
\end{webabstract}

\begin{abstract} 
For a knot $K$ in $S^3$ we construct according to Casson---or more
precisely taking into account Lin's~\cite{Lin:1992} and
Heusener's~\cite{Heu:2003} further works---a volume form on the
$\SU$--representation space of the group of $K$ (see
\fullref{Vol}). We next prove that this volume form is a
topological knot invariant (see \fullref{Invariance}) and explore
some of its properties (see \fullref{Properties}).
\end{abstract}

\maketitle

\section*{Motivation and main ideas}

In 1985, A. Casson constructed an integer valued invariant of integral
homology $3$--spheres. The original definition of Casson's invariant is
based on $\SU$--representation spaces. Informally speaking, the Casson
invariant of an homology $3$--sphere $M$ counts algebraically the
number of conjugacy classes of irreducible $\SU$--representations of
$\pi_1(M)$ in the same sense that the Lefschetz number of a map counts
the number of fixed points, see Akbulut--McCarthy~\cite{AM:1990} or
Guillou--Marin~\cite{GM:1992}.  In 1992, X-S~Lin used an analogue of
Casson's original construction to define an integer valued knot
invariant~\cite{Lin:1992}, where he indirectly proved that this
invariant is equal to half the signature of the knot. At first sight,
this equality between two apparently different quantities seems
mysterious. In 2003, M~Heusener explained Lin's result using an
orientation on the representation space of knot groups. More
precisely, given a knot $K \subset S^3$ we let $M_K$ denote its
exterior and $G_K = \pi_1(M_K)$ its group. In general the
$\SU$--representation space of $G_K$ has singularities; to avoid this
difficulty, Heusener and Klassen introduced the notion of a {regular}
representation in~\cite{HeuKlassen:1997}.  An irreducible
representation $\rho$ of $G_K$ in $\SU$ is called {regular} if
the real vector space $H^1_\rho(M_K)$ is $1$--dimensional (over
$\IR$). Here $H^*_\rho(M_K)$ denotes the $({Ad\circ \rho})$--twisted
cohomology of $M_K$. We let $\Reg{K}$ denote the set of conjugacy
classes of regular representations of $G_K$ in $\SU$. Heusener proved
that $\Reg{K}$ is a canonically \emph{oriented} $1$--dimensional
manifold (see~\cite[Section~1]{Heu:2003}).

	In this article, we investigate a {volume form} on $\Reg{K}$. We prove that it is an invariant of $K$ and we explore some of its properties. The main result of this paper is the following theorem.

\begin{maintheorem*}
	The $1$--dimensional manifold $\Reg{K}$ carries a well-defined ``ca\-no\-ni\-cal" $1$--volume form $\omega^K$. This volume form is a knot invariant. It does not depend on the orientation of $K$; and if $K^*$ denotes the miror-image of $K$, then as oriented manifolds equipped with the canonical volume form, we have $$(\Reg{K^*}, \omega^{K^*}) = (-\Reg{K}, - \omega^K).$$  
\end{maintheorem*}

	The construction which enables us to define the ``canonical" volume form $\omega^K$ on $\Reg{K}$ is motivated by the original construction of Casson's invariant. In our construction the Heegaard splitting will be replaced by a plat presentation of $K$ and we will not only compare orientations but ``natural" volume forms on some appropriate representation subspaces of $G_K$ (see~\fullref{Vol}).
	
It is possible to prove that the definition of $\omega^K$ can be
reformulated as a ``combinatorial invariant", \ie using the CW--complex
structure of the exterior of $K$. Informally speaking, this
``combinatorial invariant" is a Reidemeister torsion form on
$H^1_\rho(M_K) \cong \tangent{[\rho]}{\widehat{R}(M_K)}$. Its
definition needs Turaev's sign-determined Reidemeister torsion of
CW--complexes (see for example Turaev's monograph~\cite{Turaev:2002})
and certain distinguished bases for the twisted cohomology groups of
$M_K$. This property of $\omega^K$ is discussed by the author
in~\cite{torsionvol}. This equality between two apparently different
topological invariants---one by means of Reidemeister torsion and
another using Casson's original construction---can be considered as an
analogue of a result of E Witten about the moduli space of a
Riemannian surface. In~\cite{Witten:1991}, E Witten obtained a
remarkable formula to compute the volume of the moduli space of a
Riemannian surface in terms of a combinatorial invariant, namely a
Reidemeister torsion form on the first twisted cohomology group of the
Riemannian surface. Finally, note that this reformulation of the
volume form $\omega^K$ using the Reidemeister torsion of $M_K$ gives
another proof of its invariance.

\subsection*{Organization} 
The paper is organized as follows. Sections~\ref{RepSp}--\ref{TwistH}
consist of a review of the notions of $\SU$--repre\-sen\-ta\-tion
spaces, volume forms and regularity for a representation of a knot
group. In \fullref{Vol} we describe in detail the construction of the
``canonical" volume form on $\Reg{K}$. \fullref{Invariance} contains
the proof of the first part of the Main Theorem \ie the invariance of
the volume form (see \fullref{theoreminvariance}). In
\fullref{Properties} we finish the proof of the Main Theorem by
proving some basic properties of the volume form, we also establish a
connected sum formula. We end the paper by an explicit computation of
the volume form associated to torus knots.
	 
\section{Preliminaries}\label{RepSp}

	In this section, we collect some well-known results about
	$\SU$--representation spaces, volume forms and introduce the
	notation used throughout this paper.

\subsection{Some notation}
%	In all this article ``\emph{connected}" means ``\emph{arc-connected}".
	The fundamental group $\pi_1(W)$ of a connected CW--complex $W$ is consider without specifying a base point since all the constructions we perform are invariant under conjugation (see Porti \cite[page 9]{Porti:1997}).
	
	The Lie group $\SU$ acts on its Lie algebra $\su$ via the adjoint representation $Ad_A \co \su \to \su$ defined by $Ad_A(x) = AxA^{-1}$, where $A \in \SU$. As a manifold $\SU$ is identified with the $3$--sphere $S^3$. Furthermore we identify the $2$--sphere $S^2$ with the set of zero-trace matrices of $\SU$: $S^2 = \{A \in \SU \;|\; \mathrm{Tr}(A) = 0\}$. Recall that for each $A \in \SU$ there are $\theta \in [0, \pi]$ and $P \in S^2$ such that $A = \cos(\theta) + \sin(\theta)P$. Moreover the pair $(\theta, P)$ is {unique} if and only if $A \ne \pm \I$. Note that $Ad_A$ is the rotation of angle $2 \theta$ which fixes $P$. We always think of $\SO$ as the base space of the usual two-fold covering $\SU \to \SO$ given by $A \mapsto Ad_A$.
	
	The Lie algebra $\su$ is equipped with the usual scalar product defined by $\langle x, y \rangle = -1/2 \cdot \mathrm{Tr}(xy)$; and we identify $\su$ with the \emph{pure quaternions}, \ie with the quaternions of the form $q = a \ii + b \jj + c \kk$.

\subsection{Representation spaces}
	Given a finitely generated group $G$ we let $R(G) = \mathrm{Hom}(G; \SU)$ denote the space of $\SU$--representations of $G$. Observe that $R(G)$ is a topological space endowed with the {compact--open topology}. Here $G$ is assumed to have the discrete topology and $\SU$ the usual one. A representation $\rho\in R(G)$ is called \emph{abelian} (resp.\ \emph{central}) if its image $\rho(G)$ is an abelian subgroup of $\SU$ (resp.\ is contained in the center $\{\pm \I\}$ of $\SU$).  In the case of $\SU$, remark that a representation is abelian if and only if it is \emph{reducible} in the usual sense: $\rho$ is reducible if there exists a non-trivial proper subspace $U \subset \CC^2$ such that $\rho(g)(U) \subset U$, for all $g \in G$. A representation is called \emph{irreducible} if it is not abelian. We let $\widetilde{R}(G)$ (resp.\ $A(G)$, $C(G)$) denote  the subspace of irreducible (resp.\ abelian, central) representations. One can prove that $\widetilde{R}(G)$ is open in $R(G)$. % (resp.\ central ones). 

	The compact Lie group $\SU$ acts on $R(G)$ by conjugation. We write $[\rho]$ for the conjugacy class of the representation $\rho \in R(G)$ and we let $\SU(\rho)$ denote its orbit. The action by conjugation factors through $\SO=\SU/\{\pm \I\}$ as a {free} action on the open subspace $\widetilde{R}(G)$ and we set $\widehat{R}(G)=\widetilde{R}(G)/\SO$. In this way, we can think of the map $\widetilde{R}(G) \to \widehat{R}(G)$ as a principal $\SO$--bundle,  see Guillou--Marin~\cite[Section~3.A]{GM:1992}. 
	
\begin{notation*}
	For a connected CW--complex $W$, we write $R(W) = R(\pi_1(W))$, $\widetilde{R}(W) = \widetilde{R}(\pi_1(W))$, $\widehat{R}(W) = \widehat{R}(\pi_1(W))$ etc. 
\end{notation*}

\subsection{Representation space of knot groups}
	For a knot $K \subset S^3$ let $M_{K} = S^{3}\setminus N(K)$ denote its {exterior} and $G_{K} = \pi_{1}(M_{K})$ its {group}. Here $N(K)$ is an open tubular neighbourhood of $K$. Recall that $M_{K}$ is a compact $3$--dimensional manifold whose boundary consists in a single $2$--torus. The meridian $m$ of $K$ is only defined up to conjugation and if $K$ is oriented then $m$ is oriented by the convention $\lk(K,m) = +1$, where $\lk$ denotes the linking number.

	The abelianization $G_{K}/G_{K}' \cong H_1(M_K; \ZZ)$ is generated by the meridian $m$ of $K$. As a consequence, each {abelian} representation of $G_{K}$ is conjugate to one and only one of the $\varphi_\theta \co G_{K} \to \SU$ defined by $\varphi_\theta(m) = \cos(\theta)  + \sin(\theta) \ii$, with $0 \leqslant \theta \leqslant \pi$.

\subsection{Volume forms}
\label{S:Vol}

Here we collect some well-known facts about volume forms, see Milnor~\cite[Section 3]{Milnor:1961} and~\cite[Section 1]{Milnor:1962} for details and proofs.

\subsubsection{Volume forms and compatibility}
\label{SS:Vol}
	 
	Let $E$ be a $n$--dimensional real vector space. A \emph{volume form} $v$ on $E$ is a generator of the $n^{\rm th}$ exterior power $\bigwedge^n E^*$, where $E^* = \mathrm{Hom}_{\IR}(E, \IR)$ is the dual space of $E$.
	Let $E', E''$ be two real finite dimensional vector spaces and let $v', v''$ be volume forms on $E', E''$ respectively. The direct sum $E' \oplus E''$ inherits a canonical volume form denoted $v' \wedge v''$.
	
	Consider now a short exact sequence $\xymatrix@1@-.7pc{0 \ar[r] & E' \ar[r]^-i & E \ar[r]^-j & E'' \ar[r] & 0}$ of real finite dimensional vector spaces. Let $v'$, $v$ and $v''$ be volume forms on $E'$, $E$ and $E''$ respectively. Let $s$ denote a section of $j$ so that $i \oplus s \co E' \oplus E'' \to E$ is an isomorphism. We say that the previous three volume forms are \emph{compatible} with each other if $v' \wedge v'' = (i \oplus s)^*(v).$
	It is easy to verify that the notion of compatibility does not depend on the chosen section $s$. Finally, the following very useful lemma is quite clear.
	
\begin{lemma}\label{usefullemma}
Let  $\xymatrix@1@-.7pc{0 \ar[r] & E' \ar[r] & E \ar[r] & E'' \ar[r] & 0}$ be an exact sequence of real finite dimensional vector spaces. If any two of the vector spaces $E'$, $E$ and $E''$ are endowed with a volume form, then the third is endowed with a unique well-defined volume form which is compatible with the two others.
\end{lemma}

	In particular, if $v, v''$ are volume forms on $E, E''$ respectively, we will write $v' = v/v''$ the unique compatible volume form on $E'$ to indicate its dependence on $v$ and $v''$. Compatibility will be used in order to build up ``new" volume forms (see \fullref{SS:Construction}).

\subsubsection{The ``base $\wedge$ fiber" condition} 
\label{basefibre}
	A \emph{volume form} $v$ on a $n$--dimensional manifold is a nowhere vanishing differential $n$--form. In the sequel, we will make use of the  \emph{``base $\wedge$ fiber" condition}, which is the following. Given two volume forms $v$ and $w$ on the manifolds $M^m$ and $N^n$ respectively, a submersion $f \co M \to N$ and a point $y \in N$, then the subspace $f^{-1}(y) \subset M$ is a  submanifold of dimension $m-n$, the tangent space $\tangent{x}{f^{-1}(y)}$ is the kernel of $D_x f$ and we have the short exact sequence
\begin{equation*}\label{Base+Fibre}
\xymatrix@1@-.5pc{0 \ar[r] & \tangent{x}{f^{-1}(y)} \ar[r]^-i & \tangent{x}{M} \ar[r]^-{D_xf} & \tangent{y}{N} \ar[r] & 0.}
\end{equation*}
The submanifold  $f^{-1}(y)$ is endowed with the unique volume form $\omega$ such that, for each $x \in f^{-1}(y)$ one has $\omega_x = v_x/w_y$ \ie $\omega_x \wedge w_y = (i \oplus s)^*(v_x)$, $s$ being a section of $D_xf$.

\section{Notion of regularity}
\label{TwistH}
	
%From now on, suppose that $S^3$ and $K$ are \emph{oriented}, l
%Let $\rho : G_K \to \SU$ be any \emph{non boundary central} representation of $G_K$, \ie such that $\rho(\pi_1(\bord M_K)) \not\subset \{\pm \I\}$. 
	In this paper we will not consider the singular points of the semi-algebraic set $\widehat{R}(M_K)$ and only focus on the so-called regular representations. This section recalls the definition of {regularity} for $\SU$--representations of knot groups, see Boyer--Zhang~\cite{BZ:1996}, Heusener--Klassen~\cite{HeuKlassen:1997}, Porti~\cite[Definition 3.21]{Porti:1997} and Heusener~\cite[Section 1]{Heu:2003} for more details on this notion.   
	
\begin{notation*}
	For a CW--complex $W$ and a representation $\rho \co \pi_1(W) \to \SU$, let $H^*_\rho(W) = H^*(W; \su_{Ad \circ \rho})$ denote the cohomology of $W$ with coefficients in the adjoint representation $Ad \circ \rho$. This cohomology is called the \emph{$(Ad \circ \rho)$--twisted cohomology} of $W$. When $H^*_\rho(W) = 0$ we say that $\rho$ is \emph{acyclic}.
\end{notation*}

\subsection{Non-acyclicity of representation of knot groups}
\label{Nonacyclic}

	The long exact sequence in $(Ad \circ \rho)$--twisted cohomology corresponding to the pair $(M_K, \bord M_K)$ and Poincar\'e duality imply $\dim H^1_\rho(M_K) \geqslant 1$. So a $\SU$--representation of a knot group is {never} acyclic.
	 
	If $\rho \co G_K \to \SU$ is irreducible, then $H^0_\rho(M_K) =0$ (because $H^0_\rho(M_K)$ is equal to the subgroup of $\su$ consisting of elements fixed by $Ad \circ \rho(G_K)$). Moreover, we have $\dim H^2_\rho(M_K) = \dim H^1_\rho(M_K)$ because $$\sum_i (-1)^i \dim H^i_\rho(M_K) = 3\chi(M_K) = 0.$$ 

\subsection{\textbf{Regular representations}}
	An irreducible representation $\rho$ of $G_K$ is called \emph{regular} if $\dim H^1_\rho(M_K) = 1$. 
One can easily prove that this notion is invariant by conjugation. In the sequel, we let $\Reg{K}$ denote the set of conjugacy classes of regular representations of $G_K$ in $\SU$.

\begin{example}
	If $K$ denotes a torus knot or the figure eight knot, then any irreducible representation of $G_K$ in $\SU$ is regular (see~\cite[Example 1.43]{JDTHESE} or~\cite{torsionvol}).
\end{example}
		
	M Heusener and E Klassen proved~\cite[Proposition 1]{HeuKlassen:1997} that $\Reg{K}$ is a {$1$--dimensional manifold} (which may be empty). If $\rho$ is regular, then its conjugacy class $[\rho]$ is a smooth point of $\widehat{R}(M_K)$, $\dim \widehat{R}(M_K) = 1$ in a neighbourhood of $[\rho]$ and $\tangent{[\rho]}{\widehat{R}(M_K)}$ is isomorphic to $H^1_\rho(M_K)$ (see~\cite[Section 4]{torsionvol}). In \fullref{SubSpce}, we will see another formulation of the concept of regularity in terms of transversality of some appropriate representation subspaces of $\widehat{R}(M_K)$.

\section{Construction of the volume form}\label{Vol}

	 In this section, we explain in details the Casson-type construction of the ``canonical" {volume form} on $\Reg{K}$, see \cite{CRAS} and~\cite[Chapter 3]{JDTHESE}. First, let us make a digression on the construction of Casson's invariant. 
	 
	 Let $M$ be an oriented integral homology $3$--sphere. The original construction of Casson's invariant is based on $\SU$--representation spaces. More precisely, consider a Heegaard splitting $M = H_1 \cup_F H_2$ of $M$, where $H_i$ is a handlebody and $F = H_1 \cap H_2$ is a surface of genus $g$. The construction of Casson's invariant is based on the fact that such a splitting of $M$ gives rise to embeddings $\widehat{R}(M) \hookrightarrow \widehat{R}(H_i)$ and $\widehat{R}(H_i) \hookrightarrow \widehat{R}(F)$. Furthermore, $\widehat{R}(M)$ can be viewed as the intersection of the images of $\widehat{R}(H_1)$ and $\widehat{R}(H_2)$ inside $\widehat{R}(F)$. The crucial point is that the spaces $\widehat{R}(H_i)$ and $\widehat{R}(F)$ are canonically oriented manifolds. Informally speaking, the Casson invariant $\lambda(M)$ is the ``algebraic intersection number" of $\widehat{R}(H_1)$ and $\widehat{R}(H_2)$ in $\widehat{R}(F)$. The technical difficulties of the construction are to make sense of the algebraic intersection number of these proper open submanifolds and to show that it is independent of the Heegaard splitting (see~\cite{GM:1992} for details). 
	 
	 In our construction the Heegaard splitting will be replaced by a splitting of the knot exterior induced by a plat presentation of the knot and the orientations will be replaced by volume forms. Moreover, we will not only consider groups but {marked groups}. A  \emph{marked group} is a pair $(G, \mathcal{G})$ where $G$ is a finitely generated group and $\mathcal{G}$ is a fixed finite set of generators of $G$. We will associate to marked groups  appropriate representation subspaces. All this material will be discuss in the following subsections.
	
\begin{convention}
The $3$--sphere $S^3$ is assumed to be oriented.
\end{convention}%	From now on, we assume that $K$ denotes an oriented knot in the $3$--sphere $S^3$ which is also supposed to be oriented.

\subsection{Plat presentation and splitting of knot exteriors}
\label{platpresentation}
\begin{figure}[ht!]
\begin{center}
\includegraphics{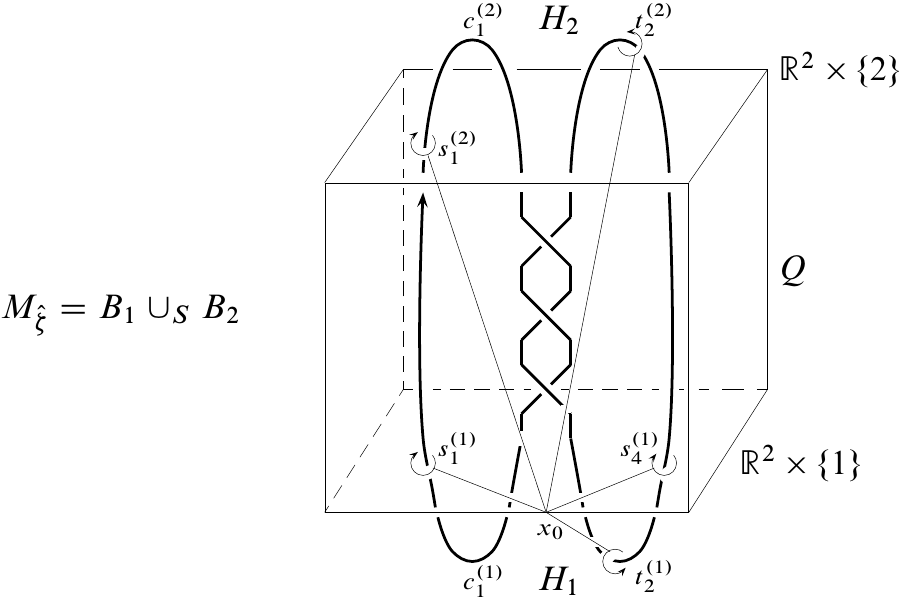}
\caption{Special systems of generators with $\varepsilon=+1$.}
\label{FigPlatGenQ}
\end{center}
\end{figure}

	Each knot $K \subset S^3$ can be presented as a $2n$--plat $\hat{\zeta}$, where $\hat{\zeta}$ is obtained from the $2n$--braid $\zeta \in B_{2n}$ by closing it with $2n$ half circles as in \fullref{FigPlatGenQ}. Explicitly, assume that the $3$--sphere $S^3 = \IR^3 \cup\{\infty\}$  is {oriented}. We choose $\epsilon \in \{\pm 1\}$ such that the basis $(e_1, e_2, e_3)$ represents the induced orientation of $\IR^3$, where $e_1 = \epsilon \ii = (\varepsilon, 0, 0)$, $e_2 = \jj = (0, 1, 0)$, $e_3 = \kk = (0, 0, 1)$. For $j=1, \ldots, 2n$, set
$$p_j = \begin{cases}
      (j,0) \in \IR^2 & \text{if $\varepsilon=1$}, \\
      (2n+1-j, 0) \in \IR^2 & \text{if $\varepsilon=-1$}.
\end{cases} $$
	Let $J = [1, 2]$,  $H_1=\{(x,y,z) \in \IR^3 \;|\; z\leqslant 1\}$, $H_2=\{(x,y,z) \in \IR^3 \;|\; z\geqslant 2\}$ and let $Q$ denote the cube $[0, 2n+1]\times [-1,1] \times J \subset \IR^2 \times J$. 
	We assume that $\zeta$ is contained in $Q$ and in a small neighbourhood of the plane $\IR \times \{0\} \times \IR$. We also assume that $\zeta \cap (\IR^2 \times \{i\}) = \mathbf{p} \times \{i\}$, $i=1,2$, where $\mathbf{p}=(p_1, \ldots, p_{2n})$. The $2n$--plat $\hat{\zeta}$ is obtained from $\zeta$ by closing it with two systems of $n$ half-circles $C_i = \{c^{(i)}_k\}_{1\leqslant k \leqslant n}\subset H_i \cap (\IR \times \{0\} \times \IR)$, $i = 1, 2$. We assume that the endpoints of $c^{(i)}_k$ are exactly $p_{2k-1}\times \{i\}, p_{2k} \times \{i\}$ in $\bord H_i$, see \fullref{FigPlatGenQ}.
	
	Such a presentation of $K$ as a $2n$--plat $\hat{\zeta}$ gives rise to a {splitting} of its exterior of the form $M_{\hat{\zeta}} = B_1 \cup_S B_2$, where $B_1, B_2$ are two handlebodies of genus $n$ and $S = B_1 \cap B_2 = S^2 \setminus N(\hat{\zeta})$ is a $2n$--punctured $2$--sphere (\cf \cite[Section 3]{Heu:2003}). To be more precise, $S= S^2 \setminus N(\hat{\zeta}) = (\IR^2 \setminus N(\mathbf{p})\times \{1\}) \cup \{\infty\}$, $B_1 = (H_1\setminus N(C_1))\cup \{\infty\}$ and $B_2=((H_2\cup \IR^2\times J )\setminus N(C_2 \cup \zeta))\cup \{\infty\}$, see \fullref{FigPlatGenQ}.

	This decomposition is similar to the Heegaard splitting used in the construction of Casson's invariant. It also gives rise to \emph{special systems of generators} (depending on the orientation of $S^3$) for $\pi_1 (B_i)$, $i=1,2$,  and $\pi_1 (S)$ respectively denoted $\mathcal{T}_i = \{t_j^{(i)}, 1\leqslant j \leqslant n\}$ and $\mathcal{S} = \{s^{(1)}_j, 1 \leqslant j\leqslant 2n\}$, see \fullref{FigPlatGenQ}. We will precisely define these generators in \fullref{choices}.

	With the notation above, the group $\pi_1(B_i)$ is the free group with basis $\mathcal{T}_i$ and the group $\pi_1(S)$ admits the finitely presentation: 
	$$\pi_1(S) = \langle s_1, \ldots, s_{2n} \; | \; s_1 \cdots s_{2n}\rangle,$$ here $s_i=s_i^{(1)}$. Furthermore, each element of $\mathcal{T}_i$ and $\mathcal{S}$ is a {meridian} of $K$. In particular all these elements are {conjugate} in $G_K$. This obvious remark will be crucial in the sequel.

	The inclusions $S \hookrightarrow B_i$ and $B_i \hookrightarrow M_K$, $i=1,2$, give rise to the following commutative diagram:  
\begin{equation}
\xymatrix@ur{
\pi_1 (S) \ar[r]^-{\kappa_1} \ar[d]_-{\kappa_2} 
							& \pi_1 (B_1) \ar[d]^-{p_1} \\
\pi_1 (B_2) \ar[r]_-{p_2} & \pi_1 (M_K)}\!\!=G_K
\label{E:diagramGroups}
\end{equation}
Each homomorphism of diagram~\ref{E:diagramGroups} is onto.
The Seifert--Van Kampen Theorem and diagram~\ref{E:diagramGroups} combine to yield the following presentation of $G_K$:
\begin{equation}\label{EQ:Presentation}
G_K=  \left\langle t^{(1)}_i, t^{(2)}_i ,\; 1\leqslant i \leqslant n \;|\;  p_1 \circ \kappa_1 (s_j) =  p_2 \circ \kappa_2 (s_j), \; 1\leqslant j \leqslant 2n-1 \right\rangle,
\end{equation}
Observe that presentation~\ref{EQ:Presentation} is a particular {Wirtinger presentation} for $G_K$.

\subsection{Choices of generators}
\label{choices}

	Here we introduce the appropriate orientation conventions and we set up the special systems of generators corresponding to a plat presentation of the knot $K$.

	Assume that the cube $Q$ is endowed with the induced orientation of the one of $\IR^3$ and choose $x_0 = (n,-1,1) \in \bord Q$ as base point. We obtain the \emph{special systems of generators} for the fundamental groups of $B_1$, $B_2$ and $S$ as follows. 
	
The generator $s^{(i)}_j$ of $\pi_1((\IR^2\setminus
N(\mathbf{p}))\times \{i\})$ is represented by a loop in $\IR^2 \times
\{i\}$ consisting of a small circle around $p_j \times \{i\}$ and the
shortest arc in $\bord Q$ connecting it to $x_0$. The circle is
oriented according to the rule: $\lk(s^{(i)}_j, L_j) = 1$, where $L_j$
denotes the oriented line $p_j \times \IR$ (the orientation points in
negative $z$--direction), see \fullref{FigPlatGenQ}. With these
choices, $$\pi_1\left((\IR^2\setminus N(\mathbf{p})\times \{i\}) \cup
\{\infty\}\right) = \left\langle s^{(i)}_1, \ldots, s^{(i)}_{2n} \;
{\displaystyle\mid}\; s^{(i)}_1 \cdots s^{(i)}_{2n}\right\rangle.$$ In order to define
the other generators we choose an orientation for the plat
$\hat{\zeta}$. We shall see later on that the construction does not
depend on this choice (see \fullref{Cor1} (3)). The generator
$t^{(i)}_k$ of $\pi_1 (H_i \setminus N(C_i))$ is represented by a loop
consisting of a small circle around $c^{(i)}_k$ and the shortest arc
in $\IR^3$ connecting it to $x_0$. The orientation of the circle is
given by the rule: $\lk(t^{(i)}_k, \hat{\zeta}) = 1$, see
\fullref{FigPlatGenQ}.
	
	Consider the homomorphism $\lambda_i \co \pi_1((\IR^2 \setminus N(\mathbf{p})) \times \{i\}) \to \pi_1(H_i \setminus N(C_i))$ induced by inclusion. We have
\begin{equation}\label{EqEps}
\lambda_i(s^{(i)}_{2k-1}) = {(t_k^{(i)})}^{\varepsilon^{(i)}_k}, \;  \lambda_i(s^{(i)}_{2k}) = {(t_k^{(i)})}^{-\varepsilon^{(i)}_k},
\end{equation}
where $\varepsilon^{(i)}_k \in \{\pm 1\}$ depends on the orientation of $\hat{\zeta}$. Observe that the $\varepsilon^{(i)}_k$ change sign simultaneously if the orientation of $\hat{\zeta}$ is reversed.

	The braid group $B_{2n}$ can be viewed as a subgroup of $\mathrm{Aut}(F_{2n})$. Here $F_{2n}$ is identified to the fundamental group $\pi_1(Q\setminus N(\zeta))$. The braid $\zeta$ induces the automorphism $\phi_\zeta \co \pi_1(Q\setminus N(\zeta)) \to \pi_1(Q\setminus N(\zeta))$ defined by $s_j^{(2)} \mapsto s_j^{(1)}$. Hence, $s_j^{(1)}$ can be viewed as a word in the generators $s_1^{(2)}, \ldots, s_{2n}^{(2)}$. One can prove (\cf diagram~\ref{E:diagramGroups}) that $$\kappa_1 \co s^{(1)}_j \mapsto \lambda_1(s^{(1)}_j) \text{ and } \kappa_2 \co s^{(1)}_j \mapsto \lambda_2 \circ \phi_\zeta(s^{(2)}_j).$$

\subsection{Representation subspaces}
\label{SubSpce}
	Corresponding to a plat presentation of $K$---which gives rise to the splitting $M_K = B_1 \cup_S B_2$ and to special systems of generators for $\pi_1(B_i)$ and $\pi_1(S)$---we introduce some \emph{special representation subspaces} for $R(B_i)$, $i=1,2$, and $R(S)$, see~\cite[Section 3]{Heu:2003}.
	
	Consider a representation $\rho \co G_K \to \SU$ and look at its restrictions $\rho_i = \rho \circ p_i$ and $\rho_S = \rho \circ p_i \circ \kappa_i$. In this way we do not obtain all the representations of the groups $\pi_1(B_i)$ or $\pi_1(S)$ because {all} the generators of these groups are conjugate to the meridian of $K$. For this reason we introduce some appropriate representation subspaces associated to the marked groups $(\pi_1(B_i), \mathcal{T}_i)$ and $(\pi_1(S), \mathcal{S})$. 
	
	Corresponding to the marked group $(G, \mathcal{G})$ we define the subset $R^{\mathcal{G}}(G)$ of $R(G) \setminus C(G)$ by setting
\begin{equation}\label{DefR}
R^{\mathcal{G}}(G)=\{\rho \in R(G)\setminus C(G) \;|\; \mathrm{Tr}({\rho(s)})=\mathrm{Tr}({\rho(t)}) \; \forall s, t \in \mathcal{G}\}.
\end{equation}
Write $\widetilde{R}^{\mathcal{G}}(G)=R^{\mathcal{G}}(G)\cap \widetilde{R}(G)$. The spaces $R^{\mathcal{G}}(G)$ and $\widetilde{R}^{\mathcal{G}}(G)$ explicitly depend on the choice of the system of generators $\mathcal{G}$. However the action by conjugation leaves $R^{\mathcal{G}}(G)$ and $\widetilde{R}^{\mathcal{G}}(G)$ invariant because the trace-function is invariant by conjugation. Thus, the quotient $\widehat{R}^{\mathcal{G}}(G)=\widetilde{R}^{\mathcal{G}}(G)/\SO$ is well-defined. Observe that $\widetilde{R}^{\mathcal{G}}(G)$ can be identified with the total space of a principal $\SO$--bundle with  $\widehat{R}^{\mathcal{G}}(G)$ as base space. 

	Here is a concrete example. Assume that $G$ is the group $G_K$ of the knot $K \subset S^3$. Assume that $\mathcal{G}$ is a finite system of generators of $G_K$ such that each element in $\mathcal{G}$ is a meridian of $K$. Then $R^\mathcal{G}(G) = R(G) \setminus C(G)$. %We also notice that $\widetilde{R}^{\mathcal{G}}(G)$ can be identify with the total space of a principal $\SO$--bundle with base space $\widehat{R}^{\mathcal{G}}(G)$. 
	
	Let $(G, \mathcal{G})$ and $(G', \mathcal{G}')$ be two marked groups. A homomorphism $\phi \co G \to G'$ is called \emph{compatible} with $\mathcal{G}$ and $\mathcal{G}'$ if $\phi(g)$ is conjugate to an element of $\mathcal{G}' \cup (\mathcal{G}')^{-1}$ for all $g \in \mathcal{G}$. If $\phi \co G \to G'$ is compatible with $\mathcal{G}$ and $\mathcal{G}'$, then $\phi$ induces a transformation $\widehat{\phi} \co \widehat{R}^\mathcal{G}(G) \to \widehat{R}^{\mathcal{G}'}(G')$.

\begin{notation*}
We write $\widehat{R}^{\mathcal{T}_i}(B_i) = \widehat{R}^{\mathcal{T}_i}(\pi_1(B_i))$, $i = 1, 2$, $\widehat{R}^{\mathcal{S}}(S) = \widehat{R}^{\mathcal{S}}(\pi_1(S))$ etc. 
\end{notation*}	

 All epimorphisms in diagram~\ref{E:diagramGroups} are compatible with the special systems of generators $\mathcal{T}_i$, $i=1, 2$, and $\mathcal{S}$ described above (because all the elements in $\mathcal{T}_i$, $i=1, 2$, and in $\mathcal{S}$ are conjugate to the meridian). Thus corresponding to diagram~\ref{E:diagramGroups} we obtain the commutative diagram: 
\begin{equation}
\xymatrix@ur{
\widehat{R}^\mathcal{S}(S) 
	&\widehat{R}^{\mathcal{T}_1}(B_1) \ar[l]_-{\widehat{\kappa}_1}\\
\widehat{R}^{\mathcal{T}_2}(B_2) \ar[u]^-{\widehat{\kappa}_2}
		& \widehat{R}(M_K)  \ar[u]_-{\widehat{p}_1} \ar[l]^-{\widehat{p}_2}}
\label{E:diagramRepsQ}
\end{equation}
All arrows in diagram~\ref{E:diagramRepsQ} are inclusions. Therefore, we can see $\widehat{R}(M_K)$ as the intersection of the images of $\widehat{R}^{\mathcal{T}_1}(B_1)$ and $\widehat{R}^{\mathcal{T}_2}(B_2)$ inside $\widehat{R}^{\mathcal{S}}(S)$. The commutative diagram~\ref{E:diagramRepsQ} is the main ingredient to define (generically) the volume form on the $\SU$--representation space of $G_K$. 

	M Heusener~\cite{Heu:2003} proved that  $\widehat{R}^{\mathcal{T}_i}(B_i)$ is a $(2n-2)$--dimensional manifold and $\widehat{R}^{\mathcal{S}}(S)$ is a $(4n-5)$--dimensional manifold. We furthermore prove in the following proposition that they also carry ``natural" volume forms. Here, ``natural" means that the volume forms are deduced from the usual ones on $\SU$ and on $(-2, 2)$.

\begin{prop}\label{VolFormRep}
	The manifold $\widehat{R}^{\mathcal{T}_i}(B_i)$ (resp.\ $\widehat{R}^{\mathcal{S}}(S)$) carries a ``natural" $(2n-2)$--volume form denoted $v^{\widehat{R}^{\mathcal{T}_i}(B_i)}$ (resp.\ a $(4n-5)$--volume form denoted $v^{\widehat{R}^{\mathcal{S}}(S)}$).
\end{prop}

\begin{proof}
We explicitly describe the volume forms on the manifolds $\widehat{R}^{\mathcal{T}_i}(B_i)$ and $\widehat{R}^{\mathcal{S}}(S)$. Their constructions are based on \fullref{usefullemma} and require the following steps.
\begin{itemize}
  \item The Lie group $\SU$ is endowed with the $3$--volume form $\eta$ induced by the basis $\{\ii, \jj, \kk\}$. Similarly we let $\eta$ denote the $3$--volume form on $\SO = \SU/\{\pm \I\}$ deduced from the one of $\SU$. Considering the fact that the trace-function $\mathrm{Tr} \co \SU\setminus\{\pm \I\} \to (-2,2)$ is a submersion, we define a $2$--volume form $\nu$ on the $2$--sphere $S^2=\{A \in \SU \;|\; \mathrm{Tr}(A)=0\}$ exploiting the ``base $\wedge$ fiber" condition. Explicitly, for each $A \in S^2$ we have the short exact sequence
\[
\xymatrix@1@-.6pc{0 \ar[r] &\tangent{A}{S^2} \ar[r] & \tangent{A}{\SU} \ar[r] & \tangent{0}{(-2,2)} \ar[r] & 0.}
\] 
In this sequence, $\SU$ is endowed with the $3$--volume form induced by $\{\ii, \jj, \kk\}$ and $(-2, 2)$ is endowed with the usual $1$--volume form. The $2$--volume form $\nu$ on $S^2$ is the unique compatible volume form with the two others. 
  \item We define a natural $(2n+1)$--volume form on $R^{\mathcal{T}_i}(B_i)$ as follows. Using the system of generators $\mathcal{T}_i$, the map $R(B_i) \to \SU^n$ defined by $\rho \mapsto \left( {\rho(t^{(i)}_1), \ldots, \rho(t^{(i)}_n)} \right)$ is an isomorphism. This isomorphism allows us to identify $R(B_i)$ with $\SU^n$. Using the natural inclusion
$ \iota \co (-2, 2) \times (S^2)^n   \to   \SU^n$ given by $$\iota(2\cos(\theta), P_1, \ldots, P_n)  =  (\cos(\theta) + \sin(\theta)P_i)_{1\leqslant i\leqslant n},$$ we identify $R^{\mathcal{T}_i}(B_i)$ with the product $(-2, 2) \times (S^2)^n$. As a consequence, the $(2n+1)$--dimensional manifold $R^{\mathcal{T}_i}(B_i)$ is endowed with a natural $(2n+1)$--volume form $v^{R^{\mathcal{T}_i}(B_i)}$, namely the one induced by the product volume form on $(-2, 2) \times (S^2)^n$.

\item We define a natural $(4n-2)$--volume form on $\widetilde{R}^{\mathcal{S}}(S)$ as follows. Let $D^*$ be the $2n$--punctured disk $S \setminus \{\infty\}$. The fundamental group of $D^*$ is the free group of rank $2n$ with basis $\mathcal{S} = \{s_1, \ldots, s_{2n}\}$ (see \fullref{FigPlatGenQ}). Let  $U$ be the subgroup normally generated by the product $s_1 \cdots s_{2n}$. We have $\pi_1 (S) = \pi_1 (D^*) /U$. Consider the map $\varphi \co R^{\mathcal{S}}(D^*) \to \SU$ defined by $\varphi(\rho) = \rho(s_1\cdots s_{2n})$ and observe that $R^{\mathcal{S}}(S)=\varphi^{-1}(\I)$. In~\cite[Lemma 3.1]{Heu:2003}, Heusener proved that $\varphi$ is surjective and that the set of critical points of $\varphi$ coincides exactly with  the set of abelian $\SU$--representations of $\pi_1 (D^*)$.
	Thus, we have the short exact sequence
$$
\xymatrix@1@-.6pc{0 \ar[r] & \tangent{\rho}{\widetilde{R}^{\mathcal{S}}(S)} \ar[r] & \tangent{\rho}{\widetilde{R}^{\mathcal{S}}(D^*)} \ar[r] & \su \ar[r] & 0.}
$$
In this sequence, $\su$ is endowed with the $3$--volume form induced by $\{\ii, \jj, \kk\}$ and $\widetilde{R}^{\mathcal{S}}(D^*)$ is endowed with a natural $(4n+1)$--volume form (observe that $\widetilde{R}^{\mathcal{S}}(D^*)$  is an open subset of $R^{\mathcal{S}}(D^*) \cong (-2,2) \times (S^2)^{2n}$ which is endowed with the product volume form).  
	Then, $\widetilde{R}^{\mathcal{S}}(S)$ is endowed with the unique $(4n-2)$--volume form $v^{\widetilde{R}^{\mathcal{S}}(S)}$ which is compatible with the two others.
  \item Finally, if $(G, \mathcal{G})$ is one of the marked groups $(\pi_1 (B_1), \mathcal{T}_1), (\pi_1 (B_2), \mathcal{T}_2)$ or $(\pi_1 (S), \mathcal{S})$, then the map $\widetilde{R}^{\mathcal{G}}(G) \to \widehat{R}^{\mathcal{G}}(G)$ is a principal $\SO$--bundle. As a consequence, we define a volume form on $\widehat{R}^{\mathcal{G}}(G)$ exploiting the ``base $\wedge$ fiber" condition. To be more precise, we have the short exact sequence
$$
\xymatrix@1@-.6pc{0 \ar[r] & \tangent{\rho}{(\SU(\rho))} \ar[r] & \tangent{\rho}{\widetilde{R}^{\mathcal{G}}(G)} \ar[r] & \tangent{[\rho]}{\widehat{R}^{\mathcal{G}}(G)}\ar[r] & 0,}
$$
where $\SU(\rho) = \{Ad_A \circ \rho \;|\; A \in \SU \}$. Observe that $\SU(\rho) \cong \SU$ is endowed with the $3$--volume form $\eta$.  So $\widehat{R}^{\mathcal{G}}(G)$ is endowed with the unique volume form which is compatible with $\eta$ and with the natural volume form we have just constructed on $\widetilde{R}^{\mathcal{G}}(G)$.\proved
\end{itemize}
\end{proof}

\subsection{Construction of the volume form}
\label{SS:Construction}

	In this subsection, we fix a $2n$--plat presentation $\hat{\zeta}$ of the {oriented} knot $K$. 
	Having in mind the result obtained in \fullref{VolFormRep}, we are ready to complete the construction of the volume form on $\Reg{\hat{\zeta}}$.

\subsubsection{Regularity and transversality}
	Let $\rho \co G_K \to \SU$ be a non-central representation and consider $\rho_i = \rho \circ p_i$ (resp.\ $\rho_S = \rho \circ p_i \circ \kappa_i$) its restriction to $\pi_1 (B_i)$ (resp.\ $\pi_1 (S)$).
	One can prove that $\rho$  is {regular} (ie, $\dim H^1_\rho(M_K) = 1$) if and only if the images in $\widehat{R}^{\mathcal{S}}(S)$ of the manifolds $\widehat{R}^{\mathcal{T}_1}(B_1)$ and $\widehat{R}^{\mathcal{T}_2}(B_2)$ intersect {transversally} at $[\rho]$, see~\cite[Proposition 3.3]{Heu:2003}. Informally speaking, the proof of this fact is essentially based on a dimensional  argument.

\subsubsection{Definition}
	The construction of the ``natural" volume form on $\Reg{\hat{\zeta}}$ combines the previous result and \fullref{VolFormRep}. It is based on the following fact: if $\rho$ is a regular representation, then %corresponding to diagram~\ref{E:diagramRepsQ}, 
\begin{equation}\label{E:SuiteExacte}
\xymatrix@1@-.5pc{0 \ar[r] & \tangent{[\rho]}{\widehat{R}(M_K)} \ar[r]^-{D_{[\rho]}\widehat{p}} & \tangent{[\rho_1]}{\widehat{R}^{\mathcal{T}_1}(B_1)} \oplus \tangent{[\rho_2]}{\widehat{R}^{\mathcal{T}_2}(B_2)} \ar[r]^-{D_{[\rho]}\widehat{\kappa}} & \tangent{[\rho_S]}{\widehat{R}^{\mathcal{S}}(S)} \ar[r] & 0}
\end{equation}
is a {short exact sequence} (cf diagram~\ref{E:diagramRepsQ} and the fact that $\widehat{R}^{\mathcal{T}_1}(B_1)$ and $\widehat{R}^{\mathcal{T}_2}(B_2)$ intersect transversally at $[\rho]$). Using the exactness of sequence~\ref{E:SuiteExacte}, we define a $1$--volume form $\vol{\Hat{\zeta}}_{[\rho]}$ on $\tangent{[\rho]}{\widehat{R}(M_{K})}$ by setting (see the notation of \fullref{basefibre}): 
\begin{equation}\label{E:DefVol}
\vol{\Hat{\zeta}}_{[\rho]} = (-1)^n \Big( v^{{\widehat{R}^{\mathcal{T}_1}(B_1)}}_{[\rho_1]} \wedge v^{{\widehat{R}^{\mathcal{T}_2}(B_2)}}_{[\rho_2]}\Big) / v^{\widehat{R}^{\mathcal{S}}(S)}_{[\rho_S]}.
\end{equation}
	In this way we locally construct a $1$--volume form $\vol{\hat{\zeta}} \co [\rho] \mapsto \vol{\hat{\zeta}}_{[\rho]} $ on the $1$--dimensional manifold $\Reg{\hat{\zeta}}$. 
	
	As it stands $\omega^{\hat{\zeta}}$ is defined in terms of the plat decomposition of the knot and apparently depends on it. In fact, we will prove  in \fullref{Invariance} (see in particular \fullref{theoreminvariance}) that it is not the case. Further observe that the normalisation given by the sign $(-1)^n$ in formula~\ref{E:DefVol} is needed to ensure the invariance of $\omega^{\hat{\zeta}}$ (see in particular \fullref{InvStab}). We will also prove in \fullref{Properties} that $\omega^{\hat{\zeta}}$ does not depend on the orientation of $\hat{\zeta}$ (see \fullref{PropProperties}). Finally remark that the order taken for the handlebodies $B_1$ and $B_2$ in formula~\ref{E:DefVol} has no importance in the definition (because of the parity of the order of the $(2n-2)$--volume form $v^{{\widehat{R}^{\mathcal{T}_i}(B_i)}}$).

%	We close this subsection by an important remark.	
\subsubsection{An important remark}\label{Remarque}
	We obtain an other splitting of the exterior of the plat $\hat{\zeta}$ by choosing the punctured sphere $S'=(\IR^2 \setminus N(\mathbf{p})\times \{2\}) \cup \{\infty\}$ and the handlebodies 
$B_1'=((H_1\cup \IR^2\times J )\setminus N(C_1 \cup \zeta))\cup \{\infty\}, \; B_2'=(H_2\setminus N(C_2))\cup \{\infty\}$, see \fullref{FigPlatGenQ}. 
Using the notation of \fullref{choices}, the epimorphisms corresponding to this splitting are given by:
$$\kappa_1' \co s^{(1)}_j \mapsto \lambda_1 \circ \phi_\zeta^{-1}(s^{(1)}_j) \text{ and } \kappa_2' \co s^{(2)}_j \mapsto \lambda_2(s^{(2)}_j).$$
Thus $\kappa_i = \kappa_i' \circ \phi_\zeta$, for $i=1,2$. Set $\widehat{Q}'_i=\widehat{\kappa}_i'(\widehat{R}^{\mathcal{T}_i}(B_i'))$, $i=1,2$, see diagram~\ref{E:diagramRepsQ}. 

The volume form $\vol{\hat{\zeta}}$ on $\Reg{\hat{\zeta}}$ can be defined using arbitrarily one of the two splittings $M_{\hat{\zeta}} = B_1 \cup_S B_2$ or $M_{\hat{\zeta}} = B'_1 \cup_{S'} B'_2$, because $\phi_\zeta \co s^{(2)}_j \mapsto s^{(1)}_j$ induces a volume preserving diffeomorphism from the regular part of $\widehat{Q}_1\cap\widehat{Q}_2$ to the one of $\widehat{Q}'_1 \cap \widehat{Q}'_2$. This result will be proved in the following subsection, see in particular \fullref{Cor1}.

\subsection{Dependence of the volume forms on the generator systems} 
	In this subsection, we analyse the dependence of the volume forms $v^{\widehat{R}^{\mathcal{T}_i}(B_i)}$ and  $v^{\widehat{R}^{\mathcal{S}}(S)}$ in terms of the generator systems $\mathcal{T}_i$ and $\mathcal{S}$ respectively.
	
	Let $F_{2n}$ denote the free group of rank $2n$ with basis $\mathcal{S}$. A braid $\zeta \in B_{2n}$ induces an automorphism $\phi_{\zeta} \co F_{2n} \to F_{2n}$, given by $\phi_{\zeta} \co s_i \mapsto g_i s_{\pi(i)}g_i^{-1}$, where $g_i \in F_{2n}$ and $\pi \in \mathfrak{S}_{2n}$ is a permutation. Observe that $\prod_{i=1}^{2n} \phi_{\zeta}(s_i) = \prod_{i=1}^{2n} s_i$. As a consequence, $\phi_{\zeta}$ is compatible with the system $\mathcal{S}$ and thus induces a diffeomorphism $\widehat{\phi_{\zeta}} \co \widehat{R}^{\mathcal{S}}(S) \to \widehat{R}^{\mathcal{S}}(S)$.
	
\subsubsection{A ``transfer lemma"}
	We begin our study by establishing a useful technical lemma. 
	
	Let $F_m  = \langle S_1, \ldots, S_m \; | \; - \, \rangle$ be a free group. Each endomorphism $\phi \co F_m \to F_m$ induces a transformation $\phi^\sharp \co R(F_m) \to R(F_m)$ and a map $\phi^{\mathrm{ab}} \co H^1(F_m; \ZZ) \to H^1(F_m; \ZZ)$, where $H^1(F_m; \ZZ) \cong \ZZ^m$. In fact, using the identification $R(F_m) \cong \SU^m$ induced by the presentation $\langle S_1, \ldots, S_m \; | \; - \, \rangle$ of $F_m$, $\phi^\sharp$ is explicitly given by
\[
(\rho(S_1), \ldots, \rho(S_m)) \mapsto (\rho \circ \phi(S_1), \ldots, \rho \circ \phi(S_m)).
\]
	For each $g \in F_m$ consider the evaluation map $\mathrm{ev}_g \co R(F_m) \to \SU$ given by $\mathrm{ev}_g(\rho) = \rho(g)$. If $\phi \co F_m \to F_m$ is an endomorphism, then $\mathrm{ev}_g \circ \phi^\sharp = \mathrm{ev}_{\phi(g)}$. With this notation, we have:

\begin{lemma}\label{lemmeVol}
	If $\phi \co F_m \to F_m$ is an endomorphism, then \[
\mathrm{ev}_{\phi(S_1)}^*(\eta) \wedge \ldots \wedge \mathrm{ev}_{\phi(S_m)}^*(\eta) = \det(\phi^{\mathrm{ab}}) \; \mathrm{ev}_{S_1}^*(\eta) \wedge \ldots \wedge \mathrm{ev}_{S_m}^*(\eta).
\]
\end{lemma}
\begin{proof}
	We just explain the main ideas which are similar to~\cite[Proposition 3.4]{AM:1990}. 
	
	The volume forms $\mathrm{ev}_{\phi(S_1)}^*(\eta) \wedge \ldots \wedge \mathrm{ev}_{\phi(S_m)}^*(\eta)$ and $\mathrm{ev}_{S_1}^*(\eta) \wedge \ldots \wedge \mathrm{ev}_{S_m}^*(\eta)$ are completely determined by their value at the trivial representation $\theta \co F_m \to \SU$ defined by $\theta(g) = \I$, for all $g \in F_m$ (because they are right-invariant). At the representation $\theta$, we have:
\[
(\mathrm{ev}_{\phi(S_1)}^*(\eta) \wedge \ldots \wedge \mathrm{ev}_{\phi(S_m)}^*(\eta)) (\theta) = \det(D_\theta \phi^{\sharp}) \, (\mathrm{ev}_{S_1}^*(\eta) \wedge \ldots \wedge \mathrm{ev}_{S_m}^*(\eta) )(\theta).
\]
Thus the proof reduces to the evaluation of $\det(D_\theta \phi^{\sharp})$. 

	Let $\mathrm{aug} \co \ZZ [F_m] \to \ZZ$ be the augmentation map $\mathrm{aug}(S_i) = 1$. Consider the Fox-matrix $A = {(a_{ij})}_{i,j}$ where $a_{ij} = \mathrm{aug}\left(\fox{\phi(S_i)}{S_j}\right) \in \ZZ$. With this notation we have $\phi^{\mathrm{ab}}(S_i) = \sum_{j=1}^m a_{ij} S_j$. If $\partial \co \su^m \to \su^m$ denotes the map given by $\partial(x_1, \ldots, x_m) = \left( \sum_{j=1}^m a_{ij}x_j \right)_{1 \leqslant i \leqslant m}$ then the diagram
\[
\xymatrix{\tangent{\theta(S_1)}{\SU} \times \cdots \times \tangent{\theta(S_m)}{\SU} \ar[r]^-{D_\theta \phi^\sharp} \ar[d]^-= & \tangent{\theta \circ \phi(S_1)}{\SU} \times \cdots \times \tangent{\theta \circ \phi(S_m)}{\SU} \ar[d]^-= \\
\su \times \cdots \times \su \ar[r]^-\partial & \su \times \cdots \times \su}
\]
commutes. Hence $\det (D_\theta \phi^\sharp) = \det (\partial) = \det (\phi^{\mathrm{ab}})$ as required.
\end{proof}

\subsubsection{Dependence of the volume form $v^{\widehat{R}^{\mathcal{T}_i}(B_i)}$}
	
	Let $(F_n, \mathcal{T})$ be the marked group $(\pi_1(B_i), \mathcal{T}_i)$, $i=1, 2$. Consider an automorphism $\phi \co F_n \to F_n$ and assume that there is a permutation $\pi \in \mathfrak{S}_n$ such that 
\begin{equation}\label{EqAuto}
\phi(t_j) = g_j t^{\varepsilon_j}_{\pi(j)} g_j ^{-1} \text{ where } g_j \in F_n \text{ and } \varepsilon_j \in \{\pm 1\}.
\end{equation}
It follows that $\mathcal{T}$ and $\mathcal{T}' = \phi(\mathcal{T})$ are compatible sets and that $R^{\mathcal{T}'}(F_n) = R^{\mathcal{T}}(F_n)$. As a consequence, $\phi$ induces a diffeomorphism $\widehat{\phi} \co \widehat{R}^{\mathcal{T}}(F_n) \to \widehat{R}^{\mathcal{T}}(F_n)$.
	
\begin{lemma}\label{lem}
%	Let $\phi \in \mathrm{Aut}(F_n)$ be given as in formula~\ref{EqAuto}. 
Let $N = \sharp\{ \varepsilon_j, 1 \leqslant j \leqslant n \; |\; \varepsilon_j = -1\}$. The induced diffeomorphism $\widehat{\phi} \co \widehat{R}^{\mathcal{T}}(F_n) \to \widehat{R}^{\mathcal{T}}(F_n)$ satisfies
\[
{\widehat{\phi}}^*(v^{\widehat{R}^{\mathcal{T}}(F_n)}) = (-1)^N v^{\widehat{R}^{\mathcal{T}}(F_n)}.
\]
Here $v^{\widehat{R}^{\mathcal{T}}(F_n)}$ denotes the natural volume form on $\widehat{R}^\mathcal{T}(F_n)$ constructed in \fullref{VolFormRep}.
\end{lemma}
\begin{proof}
Assume that $\rho(t_j) = \cos(\theta) + \sin(\theta) P^\rho_j$, with $0 < \theta < \pi$ and $P^\rho_j \in S^2$, for each $\rho \in \widetilde{R}^\mathcal{T}(F_n)$. Introduce the following maps:
\[
\mathrm{tr} \co R^{\mathcal{T}}(F_n) \to (-2, 2), \; \rho \mapsto 2\cos(\theta)  
\text{ and }
\mathrm{ax}_{t_j} \co R^{\mathcal{T}}(F_n) \to S^2, \; \rho \mapsto P^\rho_j.\]
Observe that $\mathrm{ax}_{t_j}$ is the composition of $\mathrm{ev}_{t_j}$ with the canonical projection onto $S^2$.
The automorphism $\phi$ induces a diffeomorphism ${\phi}^\dag \co {R}^{\mathcal{T}}(F_n) \to {R}^{\mathcal{T}}(F_n)$ such that $\mathrm{ax}_{t_j} \circ {\phi}^\dag = \mathrm{ax}_{\phi(t_j)}$.

 Further consider the following claims:
\begin{itemize}
  \item The map $(\mathrm{tr}, \mathrm{ax}_{t_1}, \ldots, \mathrm{ax}_{t_n}) \co {R}^\mathcal{T}(F_n) \to (-2, 2) \times (S^2)^n$ is a diffeomorphism. The volume form $v^{{R}^\mathcal{T}(F_n)}$ on ${R}^\mathcal{T}(F_n)$ is the pull-back of the product volume form on $(-2, 2) \times (S^2)^n$ by this isomorphism. 
    \item Consider $\mathfrak{t}_i \co F_n \to F_n,$ given by $t_i \mapsto t_i^{-1}, \; t_j \mapsto t_j$, if $j \ne i$. This transformation induces a diffeomorphism ${\mathfrak{t}_i}^\dag \co {R}^{\mathcal{T}}(F_n) \to {R}^{\mathcal{T}}(F_n)$ such that 
$${({\mathfrak{t}_i}^\dag)}^*(v^{{R}^\mathcal{T}(F_n)}) = 
- v^{{R}^\mathcal{T}(F_n)},\text{ for all }i.$$
   \item Consider $\mathfrak{t}_\tau \co F_n \to F_n,$  given by $t_j \mapsto t_{\tau(j)}$, where $\tau \in \mathfrak{S}_n$ is a transposition. The transformation $\mathfrak{t}_\tau$ induces a volume preserving diffeomorphism ${\mathfrak{t}_\tau}^\dag \co {R}^{\mathcal{T}}(F_n) \to {R}^{\mathcal{T}}(F_n)$.  
\end{itemize}
From these observations, we conclude
$${({\phi}^\dag)}^* (v^{{R}^\mathcal{T}(F_n)}) = (-1)^N v^{{R}^\mathcal{T}(F_n)}.$$
	The diffeomorphism ${\phi}^\dag$ induces a diffeomorphism $\widetilde{\phi} \co \widetilde{R}^\mathcal{T}(F_n) \to \widetilde{R}^\mathcal{T}(F_n)$ which makes the following diagram commutative:
\[
\xymatrix@-.6pc{\SO \ar[r] \ar[d]^-{=} & \widetilde{R}^\mathcal{T}(F_n) \ar[r] \ar[d]^-{\widetilde{\phi}} &  \widehat{R}^\mathcal{T}(F_n) \ar[d]^-{\widehat{\phi}} \\
\SO \ar[r] & \widetilde{R}^\mathcal{T}(F_n) \ar[r] &  \widehat{R}^\mathcal{T}(F_n)}
\]
This completes the proof.
\end{proof}

\subsubsection{Dependence of the volume form $v^{\widehat{R}^{\mathcal{S}}(S)}$}
Consider the marked groupw$(F_{2n}, \mathcal{S})$ where $F_{2n}$ is the fundamental group of the $2n$--punctured disk $S\setminus \{\infty\}$ and $\mathcal{S} = \{s_1, \ldots, s_{2n}\}$. Recall that $\pi_1(F) = F_{2n}/U$ where $U$ is the subgroup normally generated by the product $s_1 \cdots s_{2n}$. We are interested in an automorphism $\phi$  of $F_{2n}$ which satisfies:
\begin{enumerate}
  \item there is a permutation $\pi \in \mathfrak{S}_{2n}$ such that $\phi(s_j) = g_j s^{\varepsilon_j}_{\pi(j)} g_j ^{-1}$, where $g_j \in F_{2n}$ and  $\varepsilon_j \in \{\pm 1\}$.
  \item $\phi$ preserves the normal closure of $s_1 \cdots s_{2n} \in F_{2n}$, \ie $$\phi(s_1 \cdots s_{2n}) = g(s_1 \cdots s_{2n})^\varepsilon g^{-1}, \text{ where } g \in F_{2n} \text{ and } \varepsilon \in \{\pm 1\}.$$
\end{enumerate}
As before $\mathcal{S}$ and $\mathcal{S}' = \phi(\mathcal{S})$ are compatible sets, we have $R^{\mathcal{S}'}(S) = R^{\mathcal{S}}(S)$ and $\phi$ induces a diffeomorphism $\widehat{\phi} \co \widehat{R}^{\mathcal{S}}(S) \to \widehat{R}^{\mathcal{S}}(S)$.

\begin{lemma}\label{Lm}
Let $N = \sharp\{ \varepsilon_j, 1 \leqslant j \leqslant 2n \; |\; \varepsilon_j = -1\}$. If $N_\phi = N + (\varepsilon-1)/2$, then the induced diffeomorphism $\widehat{\phi} \co \widehat{R}^{\mathcal{S}}(S) \to \widehat{R}^{\mathcal{S}}(S)$ satisfies
\[
{\widehat{\phi}}^*(v^{\widehat{R}^{\mathcal{S}}(S)}) = (-1)^{N_\phi} v^{\widehat{R}^{\mathcal{S}}(S)}.
\]
Here $v^{\widehat{R}^{\mathcal{S}}(S)}$ denotes the natural volume form on $\widehat{R}^\mathcal{S}(S)$ constructed in \fullref{VolFormRep}.
\end{lemma}
\begin{proof}
	With the same notation as in the proof of \fullref{lem}, we introduce
\[
\mathrm{tr} \co R^{\mathcal{S}}(F_{2n}) \to (-2, 2), \; \rho \mapsto 2\cos(\theta) \text{ and }
\mathrm{ax}_{s_j} \co R^{\mathcal{S}}(F_{2n}) \to S^2, \; \rho \mapsto P^\rho_j.
\]
An inner automorphism of $F_{2n}$ induces the identity on $\widehat{R}^{\mathcal{S}}(S)$. Therefore we may assume for simplicity that $\phi(s_1 \cdots s_{2n}) = (s_1 \cdots s_{2n})^{\varepsilon}$ where $\varepsilon \in  \{\pm 1\}$. 

The automorphism $\phi$ induces two diffeomorphisms ${\phi}^\ddag \co \widetilde{R}^{\mathcal{S}}(F_{2n}) \to \widetilde{R}^{\mathcal{S}}(F_{2n})$ and $\widetilde{\phi} \co \widetilde{R}^{\mathcal{S}}(S) \to \widetilde{R}^{\mathcal{S}}(S)$. Moreover $\mathrm{ax}_{s_j} \circ {\phi}^\ddag = \mathrm{ax}_{\phi(s_j)}$. 

	We first analyse the action of ${\phi}^\ddag$ at the level of volume forms. As before $v^{{R}^{\mathcal{S}}(F_{2n})}$ is the pull-back of the product volume form on $(-2,2) \times (S^2)^{2n}$ by the diffeomorphism $(\mathrm{tr}, \mathrm{ax}_{s_1}, \ldots, \mathrm{ax}_{s_{2n}})$. Observe that by definition there is no central representation in $R^{\mathcal{S}}(F_{2n})$. 
	Let $i^{\mathcal{S}}\co R^{\mathcal{S}}(F_{2n}) \to {(\SU\setminus\{\pm \I\})}^{2n}$ be the usual inclusion $i^{\mathcal{S}}(\rho) = (\rho(s_1), \ldots, \rho(s_{2n}))$ and let $T \co {(\SU\setminus\{\pm \I\})}^{2n} \to (-2, 2)^{2n}$ be defined by $$T(A_1, \ldots, A_{2n}) = \left(\mathrm{Tr}(A_1), \ldots, \mathrm{Tr}(A_{2n})\right).$$
Thus the diagram
\[
\xymatrix@-.6pc{ {{R}^{\mathcal{S}}(F_{2n})} \ar[r]^-{i^{\mathcal{S}}} \ar[d]^-{{\phi}^\ddag} & {(\SU\setminus\{\pm \I\})}^{2n} \ar[r]^-{T} \ar[d]^-{{\phi}^\sharp} & (-2, 2)^{2n} \ar[d]^{\overline{\phi^\mathrm{ab}}} \\
{{R}^{\mathcal{S}}(F_{2n})} \ar[r]^-{i^{\mathcal{S}}} & {(\SU\setminus\{\pm \I\})}^{2n} \ar[r]^-{T} &(-2, 2)^{2n}}
\]
commutes. If $\Delta = \{(x, \ldots, x) \; |\; x \in (-2, 2)\}$, then ${(T\circ i^{\mathcal{S}})}^{-1}(\Delta) = {R}^{\mathcal{S}}(F_{2n})$. Using \fullref{lemmeVol}, we conclude that $${({\phi}^\ddag)}^*(v^{{R}^{\mathcal{S}}(F_{2n})}) = (-1)^N v^{{R}^{\mathcal{S}}(F_{2n})}.$$
	Next, consider the maps $\varphi \co {R}^{\mathcal{S}}(F_{2n}) \to \SU$ defined by $\rho \mapsto \rho(s_1 \cdots s_{2n})$ and $\Phi \co \SU \to \SU$ defined by $A \mapsto A^\varepsilon$. These maps make the diagram
\[
\xymatrix@-.6pc{ {\widetilde{R}^{\mathcal{S}}(S)} \ar[r] \ar[d]^-{\widetilde{\phi}} & {\widetilde{R}^{\mathcal{S}}(F_{2n})} \ar[r]^-{\varphi} \ar[d]^-{{\phi}^\ddag} & \SU \ar[d]^-\Phi \\
{\widetilde{R}^{\mathcal{S}}(S)} \ar[r] & {\widetilde{R}^{\mathcal{S}}(F_{2n})} \ar[r]^-{\varphi} & \SU}
\]
commutative. From this observation, we deduce 
$${\widetilde{\phi}}^* (v^{\widetilde{R}^{\mathcal{S}}(S)}) = (-1)^{N_\phi} v^{\widetilde{R}^{\mathcal{S}}(S)}.$$
	We finish the proof in the same way as in the one of \fullref{lem}.
\end{proof}

	Some consequences of the preceding lemma are summarised in:
\begin{corollary}\label{Cor1}
Let $\zeta \in B_{2n}$ be a braid such that $\hat{\zeta}$ is a knot.
\begin{enumerate}
  \item The diffeomorphism $\widehat{\phi_\zeta} \co \widehat{R}^\mathcal{S}(S) \to \widehat{R}^\mathcal{S}(S)$ induced by the $2n$--braid $\zeta$ is volume preserving.
  \item The automorphism $\phi \co F_{2n} \to F_{2n}$ defined by $\phi(s_j) = s^{-1}_{2n-j+1}$ (for all $j$) induces a diffeomorphism $\widehat{\phi} \co \widehat{R}^\mathcal{S}(S) \to \widehat{R}^\mathcal{S}(S)$ which satisfies $${\widehat{\phi}}^*(v^{\widehat{R}^\mathcal{S}(S)}) = - v^{\widehat{R}^\mathcal{S}(S)}.$$
  \item The volume form $\omega^{\hat{\zeta}}$  does not depend on the orientation of $\hat{\zeta}$.
  \item If one changes the orientation of $S^3$, then $\omega^{\hat{\zeta}}$ changes into $- \omega^{\hat{\zeta}}$.
\end{enumerate}
\end{corollary}

\begin{proof}\begin{enumerate}
  \item The automorphism $\phi_\zeta$ induced by the braid $\zeta$ verifies $N_{\phi_\zeta} = 0$ and $\varepsilon = 1$. Thus, \fullref{Lm} implies that $\widehat{\phi_\zeta}$ is volume preserving. 
  \item The automorphism $\phi$ defined by $\phi(s_j) = s^{-1}_{2n-j+1}$ verifies $N_\phi = 2n$ and $\varepsilon = -1$. Thus, \fullref{Lm} implies that ${\widehat{\phi}}^*(v^{\widehat{R}^\mathcal{S}(S)}) = - v^{\widehat{R}^\mathcal{S}(S)}$ as required.
  \item 	If we change the orientation of $\hat{\zeta}$ then the $\varepsilon^{(i)}_k \in \{\pm 1\}$ are changing sign simultaneously. Hence the $(2n-2)$--volume form $v^{\widehat{R}^{\mathcal{T}_1}(B_1)}$ (resp.\  $v^{\widehat{R}^{\mathcal{T}_2}(B_2)}$) on  $\widehat{R}^{\mathcal{T}_1}(B_1)$ (resp.\  $\widehat{R}^{\mathcal{T}_2}(B_2)$) is changed into $- v^{\widehat{R}^{\mathcal{T}_1}(B_1)}$ (resp.\  $- v^{\widehat{R}^{\mathcal{T}_2}(B_2)}$) and the $(4n-5)$--volume form $v^{\widehat{R}^{\mathcal{S}}(S)}$ is not affected. As a consequence, $\omega^{\hat{\zeta}}$ does not change (see equation~\ref{E:DefVol}).
  
 \item If we change the orientation of $S^3$, then as in (3) the volume form $v^{\widehat{R}^{\mathcal{T}_1}(B_1)}$ (resp.\  $v^{\widehat{R}^{\mathcal{T}_2}(B_2)}$) is changed into $- v^{\widehat{R}^{\mathcal{T}_1}(B_1)}$ (resp.\  $- v^{\widehat{R}^{\mathcal{T}_2}(B_2)}$). Moreover the generator $s^{(i)}_j$ in $\pi_1(S)$ is changed into $(s^{(i)}_{2n-j+1})^{-1}$, thus the $(4n-5)$--volume form $v^{\widehat{R}^{\mathcal{S}}(S)}$ on $\widehat{R}^{\mathcal{S}}(S)$ is changed into $- v^{\widehat{R}^{\mathcal{S}}(S)}$. Hence $\omega^{\hat{\zeta}}$ is changed into $-\omega^{\hat{\zeta}}$.\proved
\end{enumerate}
\end{proof}

\section{Invariance of the volume form}
\label{Invariance}

	The aim of this section is to establish that the volume form $\omega^{\hat{\zeta}}$ does not depend on the plat presentation used to define it and is a knot invariant.
	
\begin{theorem}\label{theoreminvariance}
	Let $K \subset S^3$ be a knot and let $\zeta_i \in B_{2n_i}$ be given such that the plat $\hat{\zeta}_i$ is isotopic to $K$ for $i=1, 2$. Consider the identification $\psi_i \co \Reg{K} \to \Reg{\hat{\zeta}_i}$ associated to the splitting of $M_K$ induced by the plat presentation $\hat{\zeta}_i$ of $K$ for $i = 1, 2$. Then
\[
	\psi_1^*(\omega^{\hat{\zeta}_1}) = \psi_2^*(\omega^{\hat{\zeta}_2}).
\]
\end{theorem}

	This result allows us to consider the $1$--volume form on $\Reg{K}$ associated to the knot $K$. In the sequel, we let $\omega^K$ denote this volume form on $\Reg{K}$. The proof of \fullref{theoreminvariance} is based on a theorem of Birman--Reidemeister. This theorem describes the relations between different braids which represent the same plat. It can be considered as an analogue for plats to Markov's theorem for closed braids.

\subsection{The Birman--Reidemeister Theorem}
	To state the Birman--Reide\-meis\-ter theorem we need some more definition and notation.

\begin{figure}[ht!]
\begin{center}
\includegraphics{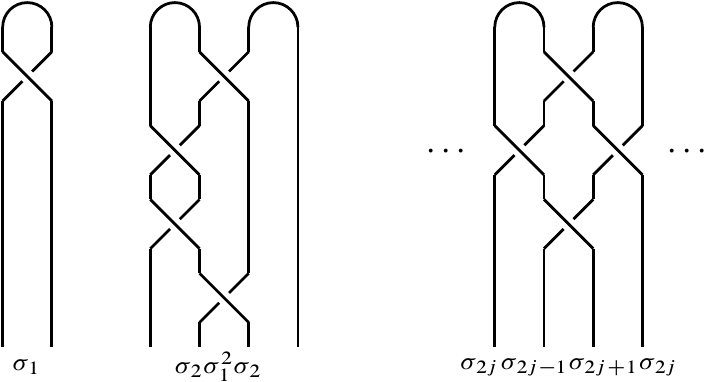}
\end{center}
\caption{The braids $\sigma_1,\; \sigma_2\sigma_1^2\sigma_2$ and $\sigma_{2j}\sigma_{2j-1}\sigma_{2j+1}\sigma_{2j}.$}\label{demitresse}
\end{figure}

	The \emph{subgroup of trivial half braid} is the subgroup $H_{2n}$ of $B_{2n}$ generated by the braids $\sigma_1, \; \sigma_2\sigma_1^2\sigma_2, \; \sigma_{2j}\sigma_{2j-1}\sigma_{2j+1}\sigma_{2j},\; 1\leq j \leq n-1$ (see \fullref{demitresse}).  It is clear that two braids of $B_{2n}$ represent the same plat if they are in the same double coset of $B_{2n}$ modulo the subgroup $H_{2n}$. 
Moreover, if $\zeta \in B_{2n}$ is such that $\hat{\zeta}$ is a knot, then $\widehat{\zeta\sigma_{2n}}$ is a knot too; and it is evident that the plats $\hat{\zeta}$ and $\widehat{\zeta\sigma_{2n}}$ are isotopic in $S^3$ (see \fullref{F:stabelem}). The operation $\zeta \mapsto \zeta \sigma_{2n}$ described  in \fullref{F:stabelem} is called an \emph{elementary stabilisation}.

	Two braids are called \emph{stably equivalent} if they represent, after a finite number of elementary stabilisations, the same double coset modulo the subgroup of trivial half braids. It is clear that two braids which are stably equivalent represent equivalent plats. The Birman--Reidemeister theorem asserts that the converse is also true:

\begin{figure}[ht!]
\begin{center}
\includegraphics{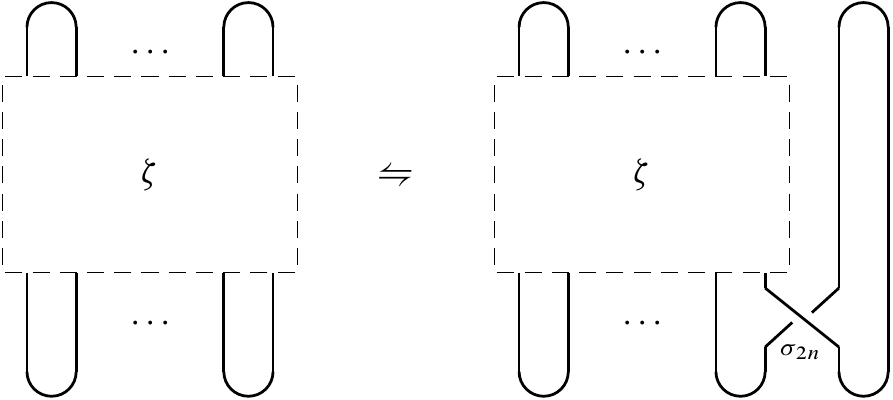}
\end{center}
\caption{Elementary stabilisation}\label{F:stabelem}
\end{figure}

\begin{theorem*}{\rm\cite{Birman:1976}}\qua
	Let $K_i \subset S^3$, $i=1, 2$, be knots and let $\zeta_i \in B_{2n_i}$ be given such that the plat $\hat{\zeta}_i$ is isotopic to $K_i$. The knots $K_1$ and $K_2$ are isotopic if and only if there exists an integer $N \geqslant \max(n_1, n_2)$ such that for each $n \geqslant N$ the braids $\zeta_i' = \zeta_i \sigma_{2n_i}\sigma_{2n_i+2} \cdots \sigma_{2n} \in B_{2n+2}$, $i=1, 2$, are in the same double coset of $B_{2n+2}$ modulo the subgroup $H_{2n+2}$.
\end{theorem*}

\subsection[Proof of \ref{theoreminvariance}]{Proof of \fullref{theoreminvariance}}
\label{proofinv}
		
	Applying the Birman--Reidemeister Theorem, the proof of \fullref{theoreminvariance} splits into two steps: we first prove the invariance of $\omega^{\hat{\zeta}}$ under the change of the double coset representative (see \fullref{InvDoublecoset}) and next prove the invariance under elementary stabilisation (see \fullref{InvStab}).
	
\subsubsection{Invariance under the change of the double coset representative}

	Write $$D_i^* = (\IR^2 \setminus N(\mathbf{p})\times \{i\}) \text{ for } i=1, 2 \text{ (see \fullref{FigPlatGenQ}).}$$  Keep the notation of \fullref{choices}. Let $F^{(i)}_{2n} = \pi_1 (D^*_i) = \langle s^{(i)}_1, \ldots, s^{(i)}_{2n} \; |\; - \,\rangle$ be the free group with basis $\mathcal{S}_i = \{ s^{(i)}_1, \ldots, s^{(i)}_{2n} \}$ and let $F^{(i)}_{n} = \pi_1 (B_i) = \langle t^{(i)}_1, \ldots, t^{(i)}_{n} \; |\; - \,\rangle$ be the one with basis $\mathcal{T}_i$, $i=1,2$, see \fullref{FigPlatGenQ}. 
	
	In~\cite{Birman:1976}, J. Birman proved that $\xi \in B_{2n}$ is in the subgroup $H_{2n}$ if and only if $\xi$ leaves the normal closure of $\{s^{(i)}_1s^{(i)}_2, \ldots, s^{(i)}_{2n-1}s^{(i)}_{2n}\}$  invariant in $F^{(i)}_{2n}$. Therefore we have an automorphism $\phi_\xi^{(i)} \co F^{(i)}_{n} \to F^{(i)}_{n}$ such that the following diagram commutes (\cf equation~\ref{EqEps}):
\[
\xymatrix{F^{(i)}_{2n} \ar[r]^-{\lambda_i} \ar[d]^-{\phi_\xi} & F^{(i)}_n \ar[d]^-{\phi_\xi^{(i)}} \\
F^{(i)}_{2n} \ar[r]^-{\lambda_i} & F^{(i)}_n}
\]
For the generators of $H_{2n}$ an elementary computation gives:
\begin{equation}
\label{E:sigma1}
\phi_{\sigma_1}^{(i)}\left(t^{(i)}_1\right) = {(t^{(i)}_1)}^{-1}, \; \phi_{\sigma_1}^{(i)}\left(t^{(i)}_j\right) = t^{(i)}_j \text{ for } j \geqslant 2\,;
\end{equation}
\begin{equation}
\label{E:sigma2}
\phi_{\sigma_2\sigma_1^2\sigma_2}^{(i)}\left(t^{(i)}_1\right) = {(t^{(i)}_2)}^{\varepsilon^{(i)}_2} t^{(i)}_1 {(t^{(i)}_2)}^{-\varepsilon^{(i)}_2}, \; \phi_{\sigma_2\sigma_1^2\sigma_2}^{(i)}\left(t^{(i)}_j\right) = t^{(i)}_j \text{ for } j \geqslant 2\,;
\end{equation}
\begin{equation}
\label{E:sigma3}
\phi_{\sigma_{2k}\sigma_{2k-1}\sigma_{2k+1}\sigma_{2k}}^{(i)}\left(t^{(i)}_j\right) = t^{(i)}_{\tau_k(j)}\text{ for } 2 \leqslant j \leqslant n.
\end{equation}
Here $\tau_k \in \mathfrak{S}_n$ denotes the transposition which permutes $k$ and $k+1$.

	Let $\zeta$ be a $2n$--braid such that $\hat{\zeta}$ is a knot. If $\xi_1, \xi_2 \in H_{2n}$ then $\widehat{\xi_1\zeta\xi_2}$ is a knot too. Let  $\kappa_i \co \pi_1 (S) \to \pi_1 (B_i)$, $\lambda_i \co F_{2n}^{(i)} \to F_{n}^{(i)}$ ($i=1,2$) denote the epimorphisms induced by $\hat{\zeta}$ and let $\kappa'_i \co \pi_1 (S) \to \pi_1 (B_i)$, $\lambda_i' \co F_{2n}^{(i)} \to F_{n}^{(i)}$ ($i=1,2$) denote the ones induced by $\widehat{\xi_1\zeta\xi_2}$.

\begin{lemma}\label{InvDoublecoset}
%	Let $\zeta$ be a $2n$--braids such that $\hat{\zeta}$ is a knot. 
For all $\xi_1, \xi_2 \in H_{2n}$, the identification $\Gamma \co \Reg{\hat{\zeta}} \to \Reg{\widehat{\xi_1\zeta\xi_2}}$ associated to the plat presentation $\widehat{\xi_1\zeta\xi_2}$ of the knot $\hat{\zeta}$ is a volume preserving diffeomorphism.
\end{lemma}
\begin{proof}
	It is sufficient to prove the lemma in the case where $\xi_i$ is one of the generators of $H_{2n}$, namely: $\xi_i = \sigma_1, \xi_i =  \sigma_2\sigma_1^2\sigma_2$ or $\xi_i = \sigma_{2k}\sigma_{2k-1}\sigma_{2k+1}\sigma_{2k} \; 1\leqslant k \leqslant n-1$. First observe that $\lambda_i'$ differs from $\lambda_i$ only if $\xi_i = \sigma_1$. Next we examine the different cases.
\begin{itemize}
  \item If $\xi_i = \sigma_1$, then $\lambda_i'(s_1^{(i)}) = {\lambda_i(s_1^{(i)})}^{-1}$ and equation~\ref{E:sigma1} implies $\kappa_i'=\kappa_i$.
  \item If $\xi_i = \sigma_2\sigma_1^2\sigma_2$, then equation~\ref{E:sigma2} implies that $\phi_{\xi_i}^{(i)} \co F^{(i)}_n \to F^{(i)}_n$ induces a volume preserving diffeomorphism $\widehat{\phi}_{\xi_i}^{(i)} \co \widehat{R}^{\mathcal{T}_i}(F^{(i)}_n) \to \widehat{R}^{\mathcal{T}_i}(F^{(i)}_n)$. Moreover, $\widehat{\kappa}_i' = \widehat{\kappa}_i \circ \widehat{\phi}_{\xi_i}^{(i)}$ (\cf diagram~\ref{E:diagramRepsQ}).
  \item If $\xi_i = \sigma_{2k}\sigma_{2k-1}\sigma_{2k+1}\sigma_{2k}$, then equation~\ref{E:sigma3} implies that the automorphism $\phi_{\xi_i}^{(i)} \co F^{(i)}_n \to F^{(i)}_n$ induces a volume preserving diffeomorphism 
  $$\widehat{\phi}_{\xi_i}^{(i)} \co \widehat{R}^{\mathcal{T}_i}(F^{(i)}_n) \to \widehat{R}^{\mathcal{T}_i}(F^{(i)}_n)$$
   (because $\nu$ is a $2$--volume form). We conclude as before.
\end{itemize} 

We deduce from these facts that there is a volume preserving
diffeomorphism\break $\psi_{\xi_i} \co \widehat{R}^{\mathcal{T}_i}(B_i) \to
\widehat{R}^{\mathcal{T}_i}(B_i)$ such that $\widehat{\kappa}_i' =
\widehat{\kappa}_i \circ \psi_{\xi_i}$, for $i=1, 2$. These
observations complete the proof because the volume form on
$\widehat{R}^{\mathcal{S}}(S)$ is not affected.
\end{proof}

\subsubsection{Invariance under stabilisation}

	Let $\zeta$ be a $2n$--braid such that $\hat{\zeta}$ is a knot. We are interested in the new braid $\zeta' = \zeta \sigma_{2n} \in B_{2n+2}$. The plat $\hat{\zeta'}$ is a knot too, and this knot is isotopic to $\hat{\zeta}$ (see \fullref{F:stabelem}). The aim of this subsection is to establish the following lemma.
\begin{lemma}\label{InvStab}	
	The identification $\Psi \co \Reg{\hat{\zeta}} \to \Reg{\widehat{\zeta\sigma_{2n}}}$ associated to the plat presentation $\widehat{\zeta\sigma_{2n}}$ of the knot $\hat{\zeta}$ is a volume preserving diffeomorphism. 
\end{lemma}
\begin{proof}
	One needs some more notation. Set $\mathcal{S}' = \mathcal{S} \cup \{s_{2n+1}, s_{2n+2}\}$. Consider the epimorphism $\kappa_i' \co \pi_1(S') \to \pi_1(B_i')$ induced by $\hat{\zeta'}$ for $i = 1, 2$ and write $Q_i = \kappa_i(R^{\mathcal{T}_i}(B_i))$, $\widehat{Q}_i = \widehat{\kappa}_i(\widehat{R}^{\mathcal{T}_i}(B_i))$ and $Q'_i = \kappa'_i(R^{\mathcal{T}'_i}(B'_i))$, $\widehat{Q}'_i = \widehat{\kappa}'_i(\widehat{R}^{\mathcal{T}'_i}(B'_i))$, see \fullref{F:stabelem}.
	
	We have $\kappa_1 = \lambda_1$, $\kappa_2 = \lambda_2 \circ \zeta$ and $\kappa_1' = \lambda_1' \circ \sigma_{2n}^{-1}$, $\kappa_2'=\lambda_2' \circ \zeta$, where $\lambda_i' \co \pi_1(S') \to \pi_1(B_i)$ is given by:
\begin{align*}
\lambda_1'(s_j^{(1)}) &=  \lambda_1(s_j^{(1)}) \text{ if } 1\leqslant j \leqslant 2n, & \lambda_2'(s_j^{(2)}) &=  \lambda_2(s_j^{(2)}) \text{ if } 1\leqslant j \leqslant 2n,  \\
\lambda_1'(s_{2n+1}^{(1)}) &= {(t_{n+1}^{(1)})}^{-\varepsilon_n^{(1)}},    & \lambda_2'(s_{2n+1}^{(2)}) &= {(t_{n+1}^{(2)})}^{-\varepsilon_n^{(2)}}, \\
\lambda_1'(s_{2n+2}^{(1)}) &= {(t_{n+1}^{(1)})}^{\varepsilon_n^{(1)}},    & \lambda_2'(s_{2n+2}^{(2)}) &= {(t_{n+1}^{(2)})}^{\varepsilon_n^{(2)}}.
\end{align*}
	The map $\lambda_i$ induces $\overline{\lambda_i} \co (-2, 2) \times (S^2)^n \to R^{\mathcal{S}}(S)$ given by
\[
\overline{\lambda_i}(2 \cos(\theta), P_1, \ldots, P_{n}) = \left({2 \cos(\theta), \varepsilon_1^{(i)} P_1, -\varepsilon_1^{(i)} P_1, \ldots, \varepsilon_n^{(i)} P_{n},-\varepsilon_n^{(i)} P_{n}}\right).
\]
Moreover $\overline{\kappa_1} = \overline{\lambda_1}$, $\overline{{\kappa_2}} = \overline{\phi_\zeta}\circ \overline{\lambda_2}$. 

We use the following notation
\[\begin{split}
(t, P_1^\zeta, \ldots, P_{2n}^\zeta)  &:= \overline{\kappa_2}(t, P_1, \ldots, P_n) \\ &\qquad = \overline{\phi_\zeta}(t, \varepsilon_1^{(2)} P_1, -\varepsilon_1^{(2)} P_1, \ldots, \varepsilon_n^{(2)} P_{n},-\varepsilon_n^{(2)} P_{n}).\end{split}
\]
We have
\begin{align*}
\overline{\kappa_1'}(2\cos(\theta), P_1, \ldots, P_{n+1}) &=
\left({2\cos(\theta), \varepsilon_1^{(1)} P_1, -\varepsilon_1^{(1)} P_1, \ldots, \varepsilon_n^{(1)} P_{n}, -\varepsilon_n^{(1)} P_{n+1},}\right. \\ & \qquad \qquad \left. {Ad_{\cos(\theta)+\varepsilon_n^{(1)}\sin(\theta)P_{n+1}}(-\varepsilon_n^{(1)}P_{n}), \varepsilon_n^{(1)}P_{n+1}}\right)
\end{align*}
and
\[
\overline{\kappa_2'}(2\cos(\theta), P_1, \ldots, P_{n+1}) =
\left({2\cos(\theta), P_1^\zeta, \ldots, P_{2n}^\zeta, -\varepsilon_n^{(1)}P_{n+1}, \varepsilon_n^{(1)}P_{n+1}}\right). 
\]

	Consider the map $f \co (-2, 2) \times (S^2)^{2n} \to  (-2, 2) \times (S^2)^{2n+2}$ given by
\begin{equation}\label{f}
f((2\cos(\theta), P_1, \ldots, P_{2n})) = (2\cos(\theta), P_1, \ldots, P_{2n}, P_{2n}, -P_{2n}).
\end{equation}
We have $f(Q_1 \cap Q_2) = Q_1' \cap Q_2'$ and $f$ induces an immersion $\widetilde{f} \co \widetilde{R}^{\mathcal{S}}(S) \to \widetilde{R}^{\mathcal{S'}}(S')$, which itself induces an immersion $\widehat{f} \co \widehat{R}^{\mathcal{S}}(S) \to \widehat{R}^{\mathcal{S'}}(S')$.

	Let $\pi_m \co (-2, 2) \times (S^2)^m \to \SU$ be the map given by
$$ \pi_m(2\cos(\theta), P_1, \ldots, P_m) = \prod_{j=1}^m (\cos(\theta) + \sin(\theta) P_j).$$ 
Observe that $\pi_{2n+2} \circ f = \pi_{2n}$. As a consequence $$D_{(t, \mathbf{P})}f_{|\pi_{2n}^*(\su)} \co \pi_{2n}^*(\su) \to \pi_{2n+2}^*(\su)$$ is a volume preserving isomorphism. Here $\pi_m^*(\su)$ denotes the pull-back.

	If $(t, \mathbf{P})=(2\cos(\theta), P_1, \ldots, P_{2n})$ corresponds to a representation in $\widetilde{R}^{\mathcal{S}}(S)$, then 
\[
\xymatrix@1{0 \ar[r] & \tangent{(t, \mathbf{P})}{\widetilde{R}^{\mathcal{S}}(S)} \ar[r]^-{D_{(t, \mathbf{P})}\widetilde{f}} & \tangent{\widetilde{f}(t, \mathbf{P})}{\widetilde{R}^{\mathcal{S'}}(S'))} \ar[r]^-{\mathrm{pr}} & \tangent{P_{2n}}{S^2} \oplus \tangent{-P_{2n}}{S^2} \ar[r] & 0}
\]
is a short exact sequence. Here $\mathrm{pr} = (\mathrm{pr}_{2n+1}, \mathrm{pr}_{2n+2})$. Hence
\begin{equation}
\label{E:T}
\tangent{\widetilde{f}(t, \mathbf{P})}{\widetilde{R}^{\mathcal{S'}}(S'))} \cong D_{(t, \mathbf{P})}\widetilde{f}(\tangent{(t, \mathbf{P})}{\widetilde{R}^{\mathcal{S}}(S)}) \oplus \tangent{P_{2n}}{S^2} \oplus \tangent{-P_{2n}}{S^2},
\end{equation}
and if $s$ denotes a section of $\mathrm{pr}$ we have, at the level of volume forms,
\[
v^{\widetilde{R}^{\mathcal{S}}(S)}_{(t, \mathbf{P})} \wedge \nu_{2n+1} \wedge \nu_{2n+2} = (D_{(t, \mathbf{P})}\widetilde{f} \oplus s)^*( v^{\widetilde{R}^{\mathcal{S}'}(S')}_{\widetilde{f}(t, \mathbf{P})}).
\]
	If $(t, \mathbf{P}) \in Q_1 \cap Q_2$ then there exist $(t, \mathbf{P}^{(i)})=(2\cos(\theta), P_1^{(i)}, \ldots, P_{n}^{(i)}) \in (-2, 2) \times (S^2)^n$ such that $\kappa_i(t, \mathbf{P}^{(i)})=(t, \mathbf{P})$, for $i=1,2$. With these preliminaries in mind, we turn to the main ingredient of the proof of \fullref{InvStab}.
\begin{claim}\label{CInvStab}
	Let $\zeta$ be a $2n$--braid such that $\hat{\zeta}$ is a knot. If $\zeta' = \zeta \sigma_{2n} \in B_{2n+2}$, then $\widehat{f} \co \widehat{R}^{\mathcal{S}}(S) \to \widehat{R}^{\mathcal{S}'}(S')$ restricts to a volume preserving diffeomorphism 
\[
\widehat{f} \co (-1)^n \widehat{Q}_1 \cap \widehat{Q}_2 \to (-1)^{n+1} \widehat{Q}'_1 \cap \widehat{Q}'_2
\]
in a neighbourhood of each regular point.
\end{claim}
\begin{proof}
	Assume that $(t, \mathbf{P}) \in Q_1 \cap Q_2$ is regular. Then $f((t, \mathbf{P}))$ is also regular in $Q'_1 \cap Q'_2$ (see \fullref{SS:Construction}). For $i = 1, 2$, 
\[
\xymatrix@1{0 \ar[r] & \tangent{(t, \mathbf{P})}{Q_i} \ar[r]^-{D_{(t, \mathbf{P})}\widetilde{f}} & \tangent{\widetilde{f}(t, \mathbf{P})}{Q_i'} \ar[r]^-{\mathrm{pr}_i} & \mathcal{V}_i \ar[r] & 0}
\]
is a short exact sequence where $\mathcal{V}_i \cong \tangent{P^{(1)}_{n}}{S^2}$.  Thus, we have
$$\tangent{\widetilde{f}(t, \mathbf{P})}{Q_i'} \cong D_{(t, \mathbf{P})}\widetilde{f}(\tangent{(t, \mathbf{P})}{Q_i}) \oplus \mathcal{V}_i.$$
Moreover, if $s_i$ denotes a section of $\mathrm{pr}_i$, then at the level of volume forms we have
\[
v^{Q_i}_{(t, \mathbf{P})} \wedge \nu_i = (D_{(t, \mathbf{P})}\widetilde{f} \oplus s_i)^*(v^{Q_i'}_{\widetilde{f}(t, \mathbf{P})}),
\]
where $\nu_i$ denotes the natural $2$--volume form on $\mathcal{V}_i$. 
	Equation~\ref{E:T} implies
$$\tangent{\widetilde{f}(t, \mathbf{P})}{\widetilde{R}^{\mathcal{S'}}(S')} \cong D_{(t, \mathbf{P})}\widetilde{f}(\tangent{(t, \mathbf{P})}{\widetilde{R}^{\mathcal{S}}(S)}) \oplus \mathcal{W}$$
where $\mathcal{W} \cong \tangent{-\varepsilon P^{(1)}_{n}}{S^2} \oplus \tangent{\varepsilon P^{(1)}_{n}}{S^2}$. As oriented vector spaces we have $-\mathcal{W} \cong \mathcal{V}_1 \oplus \mathcal{V}_2$ and $\mathcal{W}$ is endowed with the $4$--form induced by $-\nu_1 \wedge \nu_2$. Thus, the transformation $\widetilde{f} \co \widetilde{R}^{\mathcal{S}}(S) \to \widetilde{R}^{\mathcal{S'}}(S')$ induces a volume preserving diffeomorphism 
$$ \widehat{f} \co (-1)^n \widehat{Q}_1 \cap \widehat{Q}_2 \to (-1)^{n+1} \widehat{Q}'_1 \cap \widehat{Q}'_2$$
in a neighbourhood of each regular point.
\end{proof}

	This last claim achieves the proof of the invariance of $\omega^{\hat{\zeta}}$ under elementary stabilisation.
\end{proof}

\section{Some properties and applications}
\label{Properties}

\subsection{Basic properties}

	The next proposition collects some basic properties of the volume form $\omega^K$.
\begin{prop}\label{PropProperties}
	Let $K$ be an knot in $S^3$. The volume form $\omega^K$ satisfies the following assertions:
\begin{enumerate}
  \item the orientation  on $\Reg{K}$ induced by the $1$--volume form $\omega^K$ is the orientation defined in~\cite{Heu:2003},
  \item the volume form $\omega^K$ does not depend on the orientation of $K$.
  \item if $K^*$ denotes the mirror image of $K$, then the canonical isomorphism\break $\Lambda \co \Reg{K}\to \Reg{K^*}$ satisfies $\Lambda^*_{}(\omega^{K^*}) = -\omega^K$.
\end{enumerate}
\end{prop}
\begin{proof}
	Items (2) and (3) are immediate consequence of \fullref{Cor1} (3) and (4)
respectively. 
	
	Let us prove the first one. Take a $2n$--plat presentation $\hat{\zeta}$ of the knot $K$ and let $\mathfrak{O}^{\hat{\zeta}}$ denote the orientation on $\Reg{\hat{\zeta}}$ defined in~\cite{Heu:2003}. Precisely the orientation $\mathfrak{O}^{\hat{\zeta}}$ is defined~\cite[Definition 3.5]{Heu:2003} by the rule $\Reg{\hat{\zeta}} = (-1)^n \left( {\widehat{Q}_1 \cap \widehat{Q}_2} \right)$. We easily observe that the orientation induced by $v^{\widehat{R}^{\mathcal{T}_i}(B_i)}$ (resp.\ $v^{\widehat{R}^{\mathcal{S}}(S)}$) on $\widehat{R}^{\mathcal{T}_i}(B_i)$ (resp.\ $\widehat{R}^{\mathcal{S}}(S)$) corresponds to the one defined in~\cite[Sections 3 \& 4.1]{Heu:2003} (resp.\ in~\cite[Section 4.1]{Heu:2003}). Formula~\ref{E:DefVol} implies that $\mathfrak{O}^{\hat{\zeta}}$ is precisely the orientation induced by $\omega^{\hat{\zeta}}$.
\end{proof}

\subsection{A connected sum formula}
	The aim of this subsection is to prove a connected sum formula to compute the volume form $\omega^{K_1 \sharp K_2}$ associated to a composite knot ${K_1 \sharp K_2}$ in terms of the volume forms $\omega^{K_1}$ and $\omega^{K_2}$. The beginning of our study consists in a technical lemma in which we describe the regular part of the representation space of the group of $K_1 \sharp K_2$ in terms of the regular parts of the groups of $K_1$ and $K_2$.
	
	Recall that an abelian representation of $G_K$ is conjugate to one and only one of the $\varphi_\theta \co G_K \to \SU$ given by $\varphi_\theta (m) = \cos(\theta) + \sin(\theta) \ii$ with $0 \leqslant \theta \leqslant \pi$. The abelian representation $\varphi_\theta$ is called \emph{regular} if $e^{2 \theta i}$ is not a zero of the Alexander polynomial of $K$. For such representation, E. Klassen proved in~\cite[Theorem 19]{Klassen:1991} that $H^0_{\varphi_\theta}(M_K) \cong H^0(M_K; \IR) \cong \IR$ and $H^1_{\varphi_\theta}(M_K) \cong H^1(M_K; \IR) \cong \IR$. 
	
	Let $K=K_1 \sharp K_2$ be a composite knot. Let $m_i$ (resp.\ $G_i$) denote the meridian  (resp.\ the group) of $K_i$, for $i=1, 2$. Let $m$ denote the meridian of $K$. The group $G_K$ of $K$ is the amalgamated product $G_1 \ast_U G_2$. Here the subgroup $U$ of amalgamation is the normal closure of $m_1 m_2^{-1}$ in the free product $G_1 \ast G_2$ (see~\cite[Proposition 7.10]{BZ:1985}). Observe that in $G_K$ we have $m = m_1 = m_2$. 
	
	To each representation $\rho \in R(G_K)$ corresponds the restrictions $\rho_i = \rho_{|G_i} \in R(G_i)$ for $i=1, 2$. We write $\rho = \rho_1 \ast \rho_2$. If $\rho_i \in R(G_i)$ then we can form the ``composite" representation $\rho = \rho_1 \ast \rho_2$ of $G_K$ if and only if $\rho_1(m_1) = \rho_2(m_2)$. In particular, to each representation $\rho_1 \in R(G_1)$ (resp.\ $\rho_2 \in R(G_2) $) corresponds, up to conjugation, a unique abelian representation $\alpha_2 \co G_2 \to \SU$ (resp.\ $\alpha_1 \co G_1 \to \SU$) such that $\rho_1(m_1) = \alpha_2(m_2)$ (resp.\ $\alpha_1(m_1) = \rho_2(m_2)$). Hence, we can consider the following two maps:
\begin{align*}
	\iota_1 &\co \widehat{R}(M_{K_1}) \to \widehat{R}(M_K), \; [\rho_1] \mapsto [\rho_1 \ast \alpha_2],\\
	\iota_2 &\co \widehat{R}(M_{K_2}) \to \widehat{R}(M_K), \; [\rho_2] \mapsto [\alpha_1 \ast \rho_2].
\end{align*}
Using this notation, we prove:
\begin{lemma}\label{L:RegCompose} 
	Let $K=K_1 \sharp K_2$ be the connected sum of the knots $K_1$ and $K_2$. If $R_1$ (resp.\ $R_2$) denotes the space of conjugacy classes of regular representations $\rho_1 \co G_1 \to \SU$ (resp.\ $\rho_2 \co G_2 \to \SU$) such that the associated abelian representation $\alpha_2 \co G_2 \to \SU$ (resp.\ $\alpha_1 \co G_1 \to \SU$) is also regular, then:
\begin{enumerate}
  \item $R_1$ (resp.\ $R_2$) is an open submanifold of $\Reg{K_1}$ (resp.\ $\Reg{K_2}$),
  \item $\Reg{K} = \iota_1(R_1) \cup \iota_2(R_2)$,
  \item $\iota_1 \co R_1 \to \Reg{K}$ (resp.\ $\iota_2 \co R_2 \to \Reg{K}$) is an immersion.
\end{enumerate}	
\end{lemma}
\begin{proof}
	First observe that $R_i$ is obtained form $\Reg{K_i}$ by removing a finite number of points. Precisely the removed points coincide with the regular representations whose associated abelian representation corresponds to a zero of the Alexander polynomial of $K_{3-i}$ for $i=1, 2$. This is sufficient to prove the first assertion.
	
	We next prove the second item. In \cite[Proposition 12]{Klassen:1991}, E. Klassen gives the structure of the irreducible representation space of $G_K$. Any irreducible representation $\rho \in  \widetilde{R}(G_K)$ is conjugate to one of the following ``composite" representations: $\rho_1 \ast \alpha_2$, $\alpha_1 \ast \rho_2$ or $\rho_1 \ast \rho_2$, where $\rho_i \in \widetilde{R}(G_i)$ and $\alpha_i \in A(G_i)$, $i = 1, 2$. In \fullref{Fig:NoeudCompose} we illustrate our observations and we give a picture of the representation space of the group of $K = K_1 \sharp K_2$, where $K_1$ is the torus knot of type $(2, 5)$ and $K_2$ is the trefoil knot (ie, the torus knot of type $(2, 3)$). The proof  of the second assertion is based on the following two claims.

\begin{claim}
	If $\rho_i \in \widetilde{R}(G_i)$, then $\rho = \rho_1 \ast \rho_2$ is not regular.
\end{claim}
\begin{proof}[Proof of the claim]
	Corresponding to the representation $\rho = \rho_1 \ast \rho_2$  and to the amalgamated group $G_K = G_1 \ast_U G_2$ is the Mayer--Vietoris sequence
\[
\xymatrix@1@-.5pc{0 \ar[r] & \IR \ar[r]^-{\delta} & H^1_\rho(G_K) \ar[r]^-{\kappa^*} & H^1_{\rho_1}(G_1) \oplus H^1_{\rho_2}(G_2) \ar[r]^-{i^*} & \IR \ar[r] &}
\]
\[
\xymatrix@1@-.5pc{&&&&&&H^2_\rho(G_K) \ar[r] & H^2_{\rho_1}(G_1) \oplus H^2_{\rho_2}(G_2) \ar[r] & 0.}
\]
For $i = 1, 2$, we have $\dim H^1_{\rho_i}(G_i) \geqslant 1$ (see \fullref{Nonacyclic}). Hence $\rrk \kappa^* = \dim \ker i^* \geqslant 1$. As a consequence, $2 \leqslant 1 + \rrk \kappa^* = \dim H^1_\rho(G_K)$. Thus $\rho$ is not regular.
\end{proof}

\begin{claim}
	The representation $\rho_1 \ast \alpha_2$ (resp.\ $\alpha_1 \ast \rho_2$) is regular if and only if $\rho_1$ and $\alpha_2$ (resp.\ $\alpha_1$ and $\rho_2$) are regular.
\end{claim}
\begin{proof}[Proof of the claim]
	Corresponding to the irreducible representation $\rho = \rho_1 \ast \alpha_2$ and to the amalgamated group $G_K = G_1 \ast_U G_2$, the Mayer--Vietoris sequence in twisted cohomology reduces to
\[
\xymatrix@1@-.5pc{0 \ar[r] & \IR \ar[r]^-{\cong} & \IR \ar[r]^-{0} &  H^1_\rho(G_K) \ar[r]^-{\kappa^*_{(1)}} & H^1_{\rho_1}(G_1) \oplus H^1_{\alpha_2}(G_2) \ar[r]^-{i^*} & \IR \ar[r] &}
\]
\[
\xymatrix@1@-.5pc{&&&&&&&&& H^2_\rho(G_K) \ar[r]^-{\kappa^*_{(2)}} & H^2_{\rho_1}(G_1) \oplus H^2_{\alpha_2}(G_2) \ar[r] & 0.}
\]
Here $\kappa^*_{(2)}$ in onto. Hence $\dim H^2_\rho(G_K) \geqslant \dim H^2_{\rho_1}(G_1) + \dim H^2_{\alpha_2}(G_2).$

	Assume that $\rho_1$ is not regular or that $\alpha_2$ is not abelian regular. Thus we have $\dim H^2_{\rho_1}(G_1) \geqslant 2$ or $\dim H^2_{\alpha_2}(G_2) \geqslant 1$. Hence, in each case, $\dim H^2_\rho(G_K) \geqslant 2$, which proves that $\rho_1 \ast \alpha_2$ is not regular.
	
	Now, if $\rho_1$ is assumed to be regular and if $\alpha_2$ is assumed to be abelian regular, then $\dim H^1_{\rho_1}(G_1) = \dim H^2_{\rho_1}(G_1) = \dim H^1_{\alpha_2}(G_2) = 1$ and $H^2_{\alpha_2}(G_2) = 0$. Thus, $1 = \rrk i^* = 2 - \rrk \kappa^*_{(1)}$. Hence $\dim H^1_\rho(G_K) = \rrk \kappa_{(1)}^* = 1$. As a result, $\rho_1 \ast \alpha_2$ is regular. 
	
	The same arguments prove that $\alpha_1 \ast \rho_2$ is regular if and only if $\rho_2$ is regular and $\alpha_1$ is abelian regular.
\end{proof}

\begin{figure}[!ht]
\begin{center}
\includegraphics{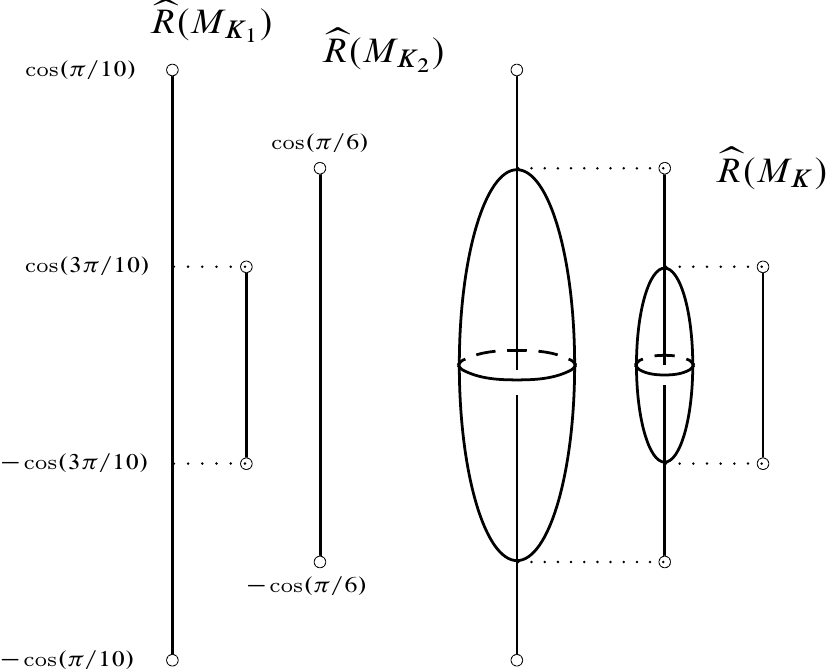}
\caption{Irreducible representations of the group of $K=K_1 \sharp K_2$.} \label{Fig:NoeudCompose}
\end{center}
\end{figure}
	Let us establish the last point. The epimorphism $\xymatrix@1@-.5pc{G_1 \ast_U G_2 \ar@{>>}[r] & G_1 \ast_U U \cong G_1}$ induces the injection $H^1_{\rho_1}(G_1) \hookrightarrow H^1_{\rho}(G_K)$. Hence, $\iota_1 \co R_1 \to \Reg{K}$ is an immersion. The same argument proves that $\iota_2 \co R_2 \to \Reg{K}$ is also an immersion.
\end{proof}

\begin{figure}[!htb]
\begin{center}
\includegraphics{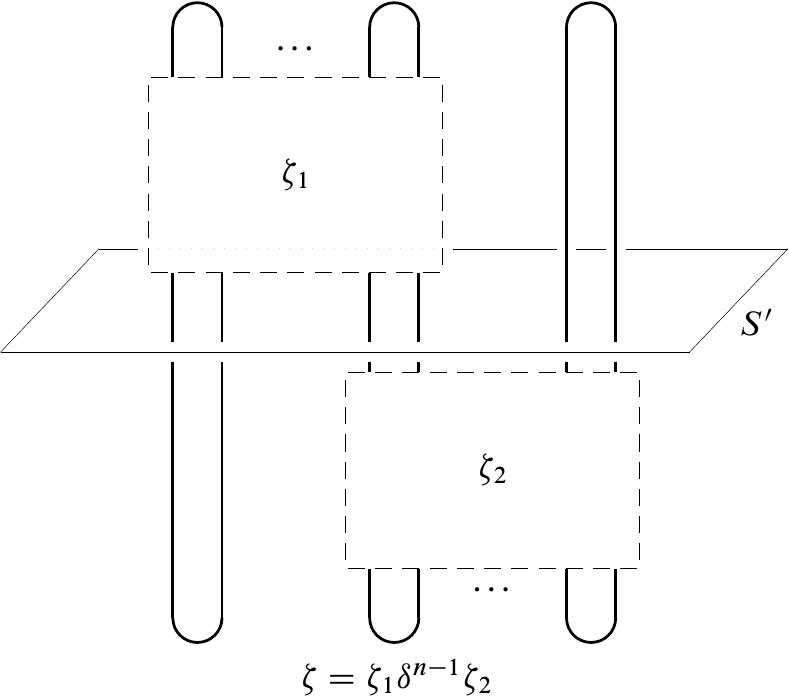}
\caption{Connected sum of plats.} \label{F:sommeconnexe}
\end{center}
\end{figure}

\begin{prop}\label{PropSom}
	Let $K=K_1 \sharp K_2$ be a composite knot. With the same notation as in \fullref{L:RegCompose}, we have
\[
\iota_1^*(\vol{K}) = \vol{K_1} \text{ and }  \iota_2^*(\vol{K}) = \vol{K_2}.
\]
\end{prop}	
\begin{proof}
	Assume that $K_1 = \hat{\zeta}_1$ is presented as a $2n$--plat and $K_2 = \hat{\zeta}_2$ as a $2m$--plat. If $\delta^{n-1}$ denotes the $(n-1)$--shift operator, then $\widehat{\zeta_1 \delta^{n-1} \zeta_2}$ is a $2(n+m-1)$--plat presentation of the connected sum $K_1 \sharp K_2$ (see \fullref{F:sommeconnexe}). \fullref{L:RegCompose} gives the form of the regular representations of $G_K$ in terms of the regular ones of $G_1$ and $G_2$. To prove $\iota_1^*(\vol{K}) = \vol{K_1}$, we use the same method as in \fullref{InvStab}. Consider the punctured $2$--sphere $\mathcal{S}' = \mathcal{S} \cup \{s_{2n+1}, \ldots, s_{2n+2(m-1)}\}$ as in \fullref{F:sommeconnexe}. This sphere allows us to split the exterior of the knot $K$. Use \fullref{Remarque} and let $$f \co (-2, 2) \times (S^2)^{2n} \to (-2, 2) \times (S^2)^{2n + 2(m-1)}$$ be the map defined by 
\begin{align*}
f((2\cos(\theta), P_1, \ldots, P_{2n})) = (2\cos(\theta), & P_1, \ldots, P_{2n}, \\
& \; \varepsilon_1 P_{2n}, -\varepsilon_1 P_{2n}, \ldots, \varepsilon_{m-1} P_{2n}, -\varepsilon_{m-1} P_{2n}).
\end{align*}
Here $\varepsilon_i \in \{\pm 1\}$ depends on the braid $\zeta_2$. As in the proof of \fullref{InvStab}, $f$ induces a transformation $\widetilde{f} \co \widetilde{R}^\mathcal{S}(S) \to \widetilde{R}^{\mathcal{S}'}(S')$. 

	The same ideas as in the proof of \fullref{InvStab} imply that $\widetilde{f}$ induces a volume preserving diffeomorphism $(-1)^n \widehat{Q}_1 \cap \widehat{Q}_2 \to (-1)^{n+m-1} \widehat{Q}_1' \cap \widehat{Q}_2'$ in a neighbourhood of each regular representation of $G_1$. We prove equality $\iota_2^*(\omega^K) = \omega^{K_2}$ by the same arguments.
\end{proof}

\subsection{An explicit computation}
	In this subsection, we give without a proof an explicit formula for the volume form $\omega^{K_{2,q}}$ associated to the torus knot $K_{2,q}$ of type $(2, q)$. Here $q\geqslant 3$ denotes an odd integer. Other explicit computations of the volume form $\omega^K$ are given in~\cite{artfibered} for fibered knots---and in particular for all torus knots---using the non abelian Reidemeister torsion. Recall that the group of the torus knot $K_{2, q}$ admits the well-known presentation $G_{K_{2, q}} = \langle x, y \; |\; x^2 = y^q\rangle$.
\begin{prop}\label{formtrefle}
Let $K_{2, q}$ be the (left-handed) torus knot of type $(2, q)$. Each irreducible representation of $G_{K_{2, q}}$ into $\SU$ is conjugate to one and only one of the $\rho_{\ell,t} \co G_{K_{2, q}} \to \SU$, where 
\begin{align*}
\rho_{\ell,t}(x) &= \ii \\
\rho_{\ell, t}(y) &= \cos\left( {(2\ell-1)\pi}/{q}\right)+\sin\left( {(2\ell-1)\pi}/{q}\right) (\cos(\pi t) \ii + \sin(\pi t)\jj),
\end{align*}
for $\ell \in \left\{1, \ldots, {(q-1)}/{2}\right\}$ and $0 < t < 1$.  

One has 
\begin{equation}
\label{Eq:Reg}
\Reg{K_{2, q}} = \bigcup_{\ell = 1}^{{(q-1)}/{2}} \left\{ {[\rho_{\ell,t}]; 0< t < 1} \right\} 
\end{equation}
and 
\begin{equation}\label{EQ:VolTrefle}
\vol{K_{2,q}}_{[\rho_{\ell, t}]}\left(\diff{\rho_{\ell, s}}{s}{s=t}\right) = \frac{8}{q} \sin^2 \left( \frac{(2\ell-1)\pi}{q} \right) \diff{\theta_m^{\rho_{\ell, s}}}{s}{s=t}
\end{equation}
where $\theta_m^{\rho_{\ell, t}} = \arccos\left((-1)^{\ell-1}\cos\left((2\ell-1)\pi/ 2q\right)\cos(\pi t) \right).$
\end{prop}

	Equation~\ref{Eq:Reg} is due to E Klassen~\cite{Klassen:1991}. One can find a proof of \fullref{formtrefle} in~\cite{torsionvol}.

\subsection{Remarks on the volume of the manifold $\Reg{K}$}

	In the introduction, we mentioned a result of Witten about the symplectic volume of the moduli space of a Riemannian surface. To close this paper we collect some remarks about the volume of $\Reg{K}$ (with respect to the canonical volume form $\omega^K$). 
	
	Because of the non-compactness of $\Reg{K}$, it is unclear if $\int_{\Reg{K}} \omega^K$ always exists in general. In the case of the torus knot of type $(2, q)$, the integral exists and we have, using \fullref{formtrefle},
\[
\int_{\Reg{K_{2,q}}} \omega^{K_{2, q}} =\sum_{\ell=1}^{\frac{q-1}{2}}(-1)^{\ell-1}
\frac{8\pi(q-2\ell+1)}{q^2} \sin^2\left( \frac{(2\ell-1)\pi}{q}\right).
\]
This formula is obtained by an elementary computation based on the knowledge of explicit descriptions of the representation space of the group of the torus knot $K_{2, q}$ and of the volume form $\omega^{K_{2, q}}$.

	Suppose that $\int_{\Reg{K}} \omega^{K}$ exists, then $\int_{\Reg{K^*}} \omega^{K^*}$ also exists and using \fullref{PropProperties} (3), we have $\int_{\Reg{K}} \omega^{K} = \int_{\Reg{K^*}} \omega^{K^*}$.

	Suppose that $\int_{\Reg{K_i}} \omega^{K_i}$ exists for $i = 1, 2$.  \fullref{PropSom} implies that $$\int_{\Reg{K_1 \sharp K_2}} \omega^{K_1 \sharp K_2}$$ also exists and that
\[
\int_{\Reg{K_1 \sharp K_2}} \omega^{K_1 \sharp K_2} = \int_{\Reg{K_1}} \omega^{K_1} + \int_{\Reg{K_2}} \omega^{K_2}.
\]

\subsection*{Acknowledgements} 
The author wishes to express his gratitude to Mi\-chael Heu\-se\-ner, Joan Porti, Daniel Lines, Vladimir Turaev, Michel Boileau, Jean-Yves LeDimet for helpful discussions related to this paper.

The author also would like to thank the referee for his careful reading of the paper and for his remarks which have contributed to making it clearer.

\bibliographystyle{gtart}
\bibliography{link}
\end{document}